\documentclass[11pt,a4paper]{article}
\usepackage{amsmath}
\usepackage{amsthm}
\usepackage{amssymb}
\usepackage{mathrsfs}
\usepackage{color,eucal}

\definecolor{ddmagenta}{rgb}{0.7,0,0.9}
\definecolor{ddcyan}{rgb}{0,0.2,1.0}
\definecolor{ddgreen}{rgb}{0,0.4,0.4}

\newcommand{\bele}{\begin{lemm}}
\newcommand{\enle}{\end{lemm}}
\newcommand{\bedef}{\begin{defi}}
\newcommand{\bete}{\begin{teor}}
\newcommand{\eddef}{\end{defi}}
\newcommand{\ente}{\end{teor}}
\newcommand{\beos}{\begin{osse}}
\newcommand{\eddos}{\end{osse}}
\newcommand{\bepr}{\begin{prop}}
\newcommand{\empr}{\end{prop}}
\newcommand{\bepro}{\begin{prob}}
\newcommand{\empro}{\end{prob}}
\newcommand{\bede}{\begin{defin}}
\newcommand{\edde}{\end{defin}}
\newcommand{\beco}{\begin{coro}}
\newcommand{\enco}{\end{coro}}

\newcommand{\beeq}[1]{\begin{equation}
 \label{#1}}
\newcommand{\eddeq}{\end{equation}}

\newcommand{\beeqa}[1]{\begin{eqnarray}
  \label{#1}}
\newcommand{\eddeqa}{\end{eqnarray}}

\newcommand{\beal}[1]{\begin{align}
 \label{#1}}
\newcommand{\eddal}{\end{align}}

\newcommand{\bespl}[1]{\begin{split}
 \label{#1}}
\newcommand{\edspl}{\end{split}}

\newcommand{\bega}[1]{\begin{gather}
 \label{#1}}
\newcommand{\edga}{\end{gather}}

\newcommand{\beeqax}{\begin{eqnarray*}}
\newcommand{\eddeqax}{\end{eqnarray*}}

\newcommand{\no}{\nonumber}

\newcommand{\tensore}{\varepsilon({\bf u})}
\newcommand{\tensoret}{\varepsilon({\mathbf{u}_t})}

\newcommand{\teta}{\vartheta}

\newcommand{\nn}{{\bf n}}
\newcommand{\uu}{{\bf u}}

\newcommand{\vv}{{\bf v}}

\newcommand{\eeta}{{\mbox{\boldmath$\eta$}}}

\newcommand{\eps}{\varepsilon}

\DeclareMathOperator{\dive}{div}

\let\TeXchi\chi
\def\chi{{\setbox0 \hbox{\mathsurround0pt
$\TeXchi$}\hbox{\raise\dp0 \copy0 }}}

\newtheorem{theorem}{Theorem}[section]
\newtheorem{corollary}{Corollary}[section]
\newtheorem{lemma}{Lemma}[section]
\newtheorem{proposition}[lemma]{Proposition}
\newtheorem{definition}[lemma]{Definition}%
\newtheorem{remark}[lemma]{Remark}%
\newtheorem{osse}[lemma]{Remark}

\newtheorem{notation}[lemma]{Notation}

\renewcommand{\part}{\partial_t}

\newcommand{\weaksto}{{\rightharpoonup^*}}
\newcommand{\weakto}{\rightharpoonup}

\newcommand{\pairing}[4]{ \sideset{_{#1 }}{_{ #2}}  {\mathop{\langle #3 , #4  \rangle}}}

 \def\fin{\hfill
         \trait .3 5 0
         \trait 5 .3 0
         \kern-5pt
         \trait 5 5 -4.7
         \trait 0.3 5 0
 \medskip}
 \def\trait #1 #2 #3 {\vrule width #1pt height #2pt depth #3pt}
\newcommand{\forae}{\text{for a.e.}}
\newcommand{\aein}{\text{a.e.\ in}}

\newcommand{\up}{\to}

\newcommand{\R}{\Bbb{R}}
\newcommand{\N}{\Bbb{N}}

\newcommand{\piecewiseConstant}[2]{\overline{#1}_{\kern-1pt#2}}

\newcommand{\underlinepiecewiseConstant}[2]{\underline{#1}_{\kern-1pt#2}}

\newcommand{\piecewiseLinear}[2]{{#1}_{\kern-1pt#2}}

\newcommand{\pwM}[2]{\widetilde{#1}_{\kern-1pt#2}}
 \def\trait #1 #2 #3 {\vrule width #1pt height #2pt depth #3pt}
\newcommand{\pwN}[2]{#1_{\kern-1pt#2}}
 \def\trait #1 #2 #3 {\vrule width #1pt height #2pt depth #3pt}

\newcommand{\dom}{\text{{\rm D}}}

\newcommand{\matrid}{\mathbf{1}}
\newcommand{\Id}{{\rm Id}}
\newcommand{\bfw}{\mathbf{W}}
\newcommand{\sig}{\sigma}
\newcommand{\traccia}{\lower3pt \hbox{$|_{\Gamma_c}$}}

\newcommand{\hunoc}{H^1 (\Gamma_c)}
\newcommand{\bfW}{\mathbf{W}}

\newcommand{\calelleveps}{\ln_\mu}
\newcommand{\ppmu}{\mathrm{p}_\mu}
\newcommand{\HH}{\mathrm{H}}
\newcommand{\HHmu}{\mathrm{H}_\mu}

\newcommand{\pepsmu}{(\mathbf{P}_{\eps}^\mu)}
\newcommand{\pepskmu}{(\mathbf{P}_{\eps_k}^\mu)}
\newcommand{\pmu}{(\mathbf{P}^\mu)}
\newcommand{\pmuk}{(\mathbf{P}^{\mu_k})}
\newcommand{\dd}{\mathrm{d}}

\newcommand{\tetaessezero}{\teta_s^0}

\newcommand{\tetame}{\teta_{\eps\mu}}
\newcommand{\tetamek}{\teta_{\eps_k\mu}}
\newcommand{\tetasme}{\teta_{s,\eps\mu}}
\newcommand{\tetasmek}{\teta_{s,\eps_k\mu}}

\newcommand{\ume}{\uu_{\eps\mu}}
\newcommand{\Thetame}{\Theta_{\eps\mu}}
\newcommand{\Thetasme}{\Theta_{s,\eps\mu}}

\newcommand{\chime}{\chi_{\eps\mu}}
\newcommand{\umek}{\uu_{\eps_k\mu}}
\newcommand{\chimek}{\chi_{\eps_k\mu}}

\newcommand{\uumu}{\uu_{\mu}}
\newcommand{\chimu}{\chi_{\mu}}
\newcommand{\uumuk}{\uu_{\mu_k}}
\newcommand{\chimuk}{\chi_{\mu_k}}

\newcommand{\tetamu}{\teta_{\mu}}
\newcommand{\tetasmu}{\teta_{s,\mu}}
\newcommand{\tetamuk}{\teta_{\mu_k}}
\newcommand{\tetasmuk}{\teta_{s,\mu_k}}

\newcommand{\tetazeroemu}{\teta_{\eps}^0}
\newcommand{\tetaessezeroemu}{\teta_{s,\eps}^0}

\newcommand{\PP}{\mathrm{\mathbf{(P)}}}
\newcommand{\Cw}{\mathrm{C}_{\mathrm{w}}^0}
\newcommand{\tetasn}{\teta_{s,n}}
\newcommand{\tetasinf}{\teta_{s,\infty}}

\newcommand{\Imu}{\mathscr{I}_\mu}
\newcommand{\cdr}{(\cdot,r)}

\newenvironment{rult}{\color{ddmagenta}}{\color{black}}

\newcommand{\brult}{\begin{rult}}
\newcommand{\eult}{\end{rult}}

\begin{document}

                                %
\title{\Large Long-time behaviour
of a thermomechanical model for adhesive
 contact}

\author{
 Elena Bonetti\\
  {\small Dipartimento di Matematica ``F.\
Casorati''}\\
  {\small Universit\`a di Pavia}\\
  {\small Via Ferrata 1,  I-- 27100 Pavia,
Italy}\\
  {\small \texttt{elena.bonetti@unipv.it}}\\
  \and
  Giovanna
Bonfanti\\
  {\small Dipartimento di Matematica}\\
  {\small Universit\`a di
 Brescia}\\
  {\small Via Valotti 9, I--25133 Brescia, Italy}\\
  {\small \texttt{bonfanti@ing.unibs.it}}\\
  \and
 Riccarda Rossi\\
  {\small Dipartimento di Matematica}\\
  {\small Universit\`a di
 Brescia}\\
  {\small Via Valotti 9, I--25133 Brescia, Italy}\\
  {\small \texttt{riccarda.rossi@ing.unibs.it}}\\
}

\date{April 03, 2009}

\maketitle

 \numberwithin{equation}{section}

\begin{abstract}
 This paper deals with the large-time  analysis  of a PDE system
modelling contact  with adhe\-sion,  in the case when thermal
effects are taken into account. The  phenomenon of adhesive contact
is described in terms of phase transitions for a surface damage
model proposed by  \textsc{M. Fr\'emond}.
 Thermal effects are governed by  entropy balance laws. The
resulting system is highly nonlinear, mainly due to the presence of
internal constraints on the phy\-sical variables and the coupling of
equations written in a domain and   on  a contact surface. We prove
existence of solutions on the whole time interval $(0,+\infty)$ by a
double approximation procedure. Hence, we  are able to  show  that
solution trajectories admit cluster points which fulfil  the
stationary problem associated with the evolutionary system, and that
in  the large-time limit dissipation vanishes.
\end{abstract}
\vskip3mm


                                %
                                %
\pagestyle{myheadings} \markright{\it Bonetti, Bonfanti, Rossi /
Analysis of a thermomechanical  model for adhesive contact}
                                %

                                %

\section{Introduction}
This paper is concerned with the large-time analysis of a PDE system
describing adhesive contact between a thermo-viscoelastic body  and
a rigid support. The model has been recently introduced and
global-in-time existence results have been proved on finite-time
intervals both in isothermal cases   (see \cite{bbr1}
 in the case of an irreversible damage evolution on the contact
 surface, and~\cite{bbr2} for the reversible case),
 and   for
PDE systems  including thermal effects (see \cite{bbr3}). The
modelling approach for contact with adhesion which we apply refers
to a damage theory described by phase transitions, and it is due to
\textsc{M. Fr\'emond} (see \cite{fre}). The idea consists in
describing the adhesion between viscoelastic bodies in terms of  a
surface damage theory, in which the damage parameter is related to
the active bonds which are responsible for the adhesion between the
bodies. Hence, the equations of the evolutionary  system are
recovered from thermomechanical laws and they are written in the
domain of the viscoelastic body and on the contact surface.

 It turns
out to be interesting, both from  a theoretical point of view and in
view of applications,  to investigate how the thermomechanical
system (i.e.,  the body and the rigid support it is in contact with)
behaves  for large times. More precisely, we shall investigate   if
the trajectories of the solutions to the resulting PDE system
present some cluster point,  in the limit as time goes to $+\infty$.
Then, we  shall  look for a relation between these limit states and
the stationary system associated with our evolution problem. In
particular, we aim to prove that, in  the limit,  solution
trajectories  reach a thermomechanical equilibrium state in which
dissipation vanishes.
 This kind of  large-time analysis was  performed in
\cite{bbr2} for the reversible  model in the isothermal case.

Before introducing the long-time behaviour analysis of the problem,
we shall briefly recall the model and  make some comments on the
existence of solutions on the whole time interval~$(0,+\infty)$.
Moreover,   we shall point out
 that this paper also presents a novelty in the
formulation of the model itself,  as we generalize
 the convex potential usually ensuring internal
constraints on the damage parameter.

\subsection*{The model}
We mainly refer to the recent contribution  \cite{bbr3}, in which
(a slightly different version of)
 the thermomechanical model has been introduced.  The state
variables, in terms of which the equilibrium of the system is
established,  are defined in the domain $\Omega\subset\R^3$ (where
  the body is located),  and on the contact surface $\Gamma_c$.
 Namely,  we shall take $\Omega$ to be a sufficiently smooth
bounded domain in~$\R^3$, with boundary
$\partial\Omega=\bar\Gamma_1\cup\bar\Gamma_2\cup\bar\Gamma_c$. The
sets $\Gamma_i$ are open subsets in the relative topology of
$\partial\Omega$, with smooth boundary and disjoint one from each
other. In particular, $\Gamma_c$ is the contact surface. We suppose
that $\Gamma_c$ and $\Gamma_1$ have strictly positive measures and,
for the sake of simplicity, we identify $\Gamma_c$ with a subset of
$\R^2$, i.e.,  we shall treat  $\Gamma_c$ as a flat surface.

The state variables we shall consider in $\Omega$  are  the absolute
temperature $\teta$ of the body and the macroscopic deformations,
given in terms of the linearized strain tensor $\varepsilon({\bf
u})$ (${\bf u}$ represents the vector of small displacements). On
the contact surface, we introduce the surface absolute temperature
(the reader may think of  the temperature of the adhesive glue)
$\teta_s$ and a damage parameter $\chi$,  related to the active
bonds in the glue ensuring adhesion. For the moment, we do not
require any constraints on the values assumed by $\chi$. Taking into
account  local interactions (in the glue and between the glue and
the body) we include  the gradient $\nabla\chi$ and the displacement
trace ${\bf u}_{|_{\Gamma_c}}$ among the state variables on the
contact surface. The free energy in $\Omega$ is written as follows
\begin{equation}\label{freeO}
\Psi_\Omega=\teta(1-\ln(\teta))+\teta\hbox{tr }\varepsilon({\bf
u})+\frac 1 2 \varepsilon({\bf u})K\varepsilon({\bf u}),
\end{equation}
where $K$ is the elasticity tensor and the coefficient $\teta$
multiplying $\mathrm{tr}\varepsilon({\bf
u})$
accounts for the  thermal
expansion energy.

\begin{remark}\label{scelgolog} \upshape
Notice that here  we have taken the term $\teta(1-\ln(\teta))$ for
the purely thermal contribution in the free energy $\Psi_\Omega$
(and, similarly,  for the free energy $\Psi_{\Gamma_c}$ below),
while in \cite{bbr3} we have considered a more general concave
function. The particular choice in this paper is very frequent in
the literature, as it has  some analytical and modelling advantages.
From the latter viewpoint, the presence of the logarithm in
\eqref{freeO} yields an internal constraint on the temperature:
indeed,  the domain of $\Psi_\Omega$ is given for $\teta>0$, which
is in agreement with thermodynamical consistency. On the analytical
level, this  form of the thermal contribution   shall allow us to
simplify the procedure exploited in \cite{bbr3} to prove existence
of solutions (see Remark \ref{piufacile}).
\end{remark}

Next,  we specify  the free energy in $\Gamma_c$, which presents
some novelty with respect to the model introduced in \cite{bbr3}
(cf. also \cite{bbr1} and \cite{bbr2}).   In fact, we shall consider
\begin{equation}\label{freeG}
\Psi_{\Gamma_c}=\teta(1-\ln(\teta_s))+\lambda(\chi)(\teta_s-\teta_{eq})+\widehat\beta(\chi)+\sigma(\chi)+\frac
1 2|\nabla\chi|^2+\frac 1 2\chi^+|{\bf u}_{|_{\Gamma_c}}|^2+I_-({\bf
u}_{|_{\Gamma_c}}\cdot{\bf n}),
\end{equation}
where $\teta_{eq}$ is a critical temperature, and $\widehat\beta$ is
a  proper, convex, and lower-semicontinuous function. The indicator
function $I_-$ forces the scalar product ${\bf
u}_{|_{\Gamma_c}}\cdot{\bf n}$ to be non-positive, as it is defined
on $\R$ by $I_-(y)=0$ if $y\leq  0$ and $I_-(y)=+\infty$ for $y>0$.
 This renders  the impenetrability condition between the body
and the support. In  the same way,  the term $\widehat\beta$ may
yield a constraint on the values assumed by $\chi$. For example, a
proper choice of $\widehat{\beta}$ may enforce positivity of $\chi$
(see \cite{bbr1, bbr2, bbr3}). In particular, this occurs when,
classically, $\widehat\beta=I_{[0,1]}$, forcing $\chi\in[0,1]$.
Then, the coefficient of $|{\bf u}_{|_{\Gamma_c}}|^2$ remains
non-negative, in accord with physical consistency.  However, in the
present paper we shall allow the potential $\widehat\beta$ to be
more general and we do not impose any a priori restriction on its
domain. Hence,
 to ensure physical consistency, we
introduce a constraint on the deformation coefficient, which is
fixed to be   $\chi^+$ (using the notation $r^+ = \max(r,0)$ for
every $r \in \R$). Finally,
 the function $\lambda$ is related to the latent
heat, while $\sigma$ takes into account  possibly
non-convex contributions in the free energy. In particular, we
include in $\sigma$ cohesive effects in the glue, which are
represented by a non-increasing function in $\chi$ (a simple choice
is $\sigma(\chi)=w(1-\chi)$, with a positive parameter $w$).

Then, the evolution of the system is governed by two convex
potentials (non-negative and assuming their minimum $0$ if there is
no dissipation), namely the pseudo-potentials of dissipation written
in $\Omega$ and in $\Gamma_c$. We have (here $K_v$ is a viscosity
matrix)
\begin{equation}\label{pseudoO}
\Phi_\Omega=\frac 1 2 |\nabla\teta|^2+\frac 1 2\varepsilon({\bf
u}_t)K_v\varepsilon({\bf u}_t),
\end{equation}
and
\begin{equation}\label{pseudoG}
\Phi_{\Gamma_c}=\frac 1 2|\nabla\teta_s|^2+\frac 1 2|\chi_t|^2+\frac
1 2k(\chi)(\teta_{|_{\Gamma_c}}-\teta_s)^2.
\end{equation}
Notice that the dissipation in $\Omega$ depends on $\varepsilon({\bf
u}_t)$ and on $\nabla \teta$, while the dissipation in $\Gamma_c$
depends on $\nabla \teta_s$, $\chi_t$, and on
$(\teta_{|_{\Gamma_c}}-\teta_s)$. The function $k$, which  accounts
for the heat exchange between the body and the adhesive material,
shall be taken non-negative and smooth enough. Actually, to
characterize the large-time behaviour of the system, we need to
assume that $k$  is bounded from below by some positive constant
(see Remark \ref{rimane}).

\subsection*{The PDE system}
Proceeding as in \cite{bbr3},  we refer to thermomechanical
laws and, after specifying the constitutive equations in terms of
the above potentials, we arrive at  the following PDE system ($T$ is
a fixed final time)
\begin{align}
\label{e1} &\partial_t (\ln(\teta)) - \dive (\uu_t) -\Delta \teta= h
\qquad \text{in $\Omega \times (0,T)$,}
\\
\label{condteta} &\partial_n \teta= \begin{cases} 0 & \text{in
$(\partial \Omega \setminus \Gamma_c) \times (0,T)$,}
\\
-k(\chi) (\teta-\teta_s) & \text{in $ \Gamma_c \times (0,T)$,}
\end{cases}
\\
\label{eqtetas} &\partial_t (\ln(\teta_s)) - \partial_t
(\lambda(\chi))  -\Delta \teta_s= k(\chi) (\teta-\teta_s) \qquad
\text{in $\Gamma_c \times (0,T)$,}
\\
\label{condtetas} &\partial_n \teta_s =0 \qquad \text{in $\partial
\Gamma_c \times (0,T)$,}
\\
\label{eqI}
&-\hbox{div }(K\tensore+ K_v\tensoret+ \teta \matrid)={\bf f}\qquad\hbox{in }\Omega\times (0,T),\\
&{\bf u}={\bf 0}\quad\hbox{in }\Gamma_{1}\times (0,T),
\quad(K\tensore+K_v\tensoret + \teta \matrid){\bf n}={\bf g}\qquad\hbox{in }\Gamma_{2}\times (0,T),\label{condIi}\\
& (K\tensore+ K_v\tensoret + \teta \matrid){\bf n}+\chi^+{\bf u}+
\partial I_-(\uu\cdot {\bf n}){\bf n}
\ni{\bf 0}\qquad\hbox{in }\Gamma_{c} \times (0,T),\label{condIii}\\
&\chi_t -\Delta\chi+\beta(\chi) +
\sigma'(\chi)-\lambda'(\chi) \teta_{\rm{eq}} + \HH(\chi)\frac 1
2\vert{\bf u}\vert^2
 \ni -\lambda'(\chi)\teta_s\qquad\hbox{in }\Gamma_c
\times (0,T),\label{eqII}
\\
&\partial_{\nn_s} \chi=0 \qquad \text{in $\partial \Gamma_c \times
(0,T)\,$} \label{bord-chi}
\end{align}
where
$h$ is an external entropy source, $\mathbf{f}$ a volume force, and $g$ a traction. Moreover,
 $\beta=\partial\widehat\beta$ and $\HH = \partial \mathrm{p}$, where $\mathrm{p}(r)=r^+$
 for all $r   \in  \R$. Hence, the Heaviside maximal monotone operator $\HH: \R \to 2^{\R}$
 is defined by $\HH(r) =0$ if $r<0$, $ \HH(0) = [0,1]$, and $\HH(r) =
 1$ if $r>0$. We warn that,
  here  and in what
follows, we shall omit for simplicity the index $v\traccia$ to
denote the trace  on $\Gamma_c$ of a function $v$, defined in
$\Omega$.

\begin{remark}
\upshape Let us comment on the above equations, while referring  the
 to \cite{bbr3} for their rigorous derivation. First of all, we  point out that \eqref{e1} and
\eqref{eqtetas} are entropy equations. The possibility of describing
thermal effects in phase transitions by the use of an entropy
equation, in place of the more standard energy balance, has only
recently been introduced. In particular, let us point out that the
entropy in $\Omega$ is defined as $\ln(\teta)-\mathrm{div} ({\bf
u})$ and in $\Gamma_c$ as $\ln(\teta_s)-\lambda(\chi)$. Using
entropy equations brings to some advantages both for the analytical
treatment and the modelling of the phenomenon. In particular,   from
\eqref{e1} and \eqref{eqtetas} one directly recovers the positivity
of the temperature, which is necessary for thermodynamical
consistency, avoiding the application of any maximum principle
argument. We do not enter the details of this theory and refer,
among the others, to the papers \cite{bcf} and \cite{BFR}. Then,
\eqref{eqI} is derived from the momentum balance, in which
accelerations are not taken into account. Equation \eqref{eqII} is
recovered as a balance equation for micro-movements related to the
evolution of the phase parameter (see \cite{fre} for the theory of
the generalized principle of virtual power including micro-movements
and micro-forces responsible for the phase transition).
\end{remark}

\begin{remark}\label{rimane}
 \upshape We emphasize  that the structure of \eqref{eqII} is more
complicate than the analogous equation in the model studied in
\cite{bbr3}. Indeed, the maximal monotone operator  $\beta$ is more
general than the one considered in~\cite{bbr3}, for we do not impose
any restriction on its domain. Moreover,  the presence of the
operator $\HH$ in \eqref{eqII} introduces a new nonlinearity in the
equation.
  From the physical viewpoint,
since $\HH(0)=[0,1]$, a residual influence of macroscopic
displacements on the mechanical behaviour of the glue may persist
also when the glue is damaged, e.g. if $\chi=0$ ($\chi$ representing
here the proportion of active bonds). This corresponds to assuming
that a local interaction between the body and the support is
preserved even when the bonds in the glue are completely damaged.
Notice that this  is reasonable,   if one takes into account
distance forces. An analogous argument justifies the assumption that
$k$ in \eqref{eqtetas}  is bounded from below by some positive
constant, see~\eqref{hyp-k} later on.  This ensures that the thermal
local interaction between the body and the support  is conserved
even when the adhesion is not active.
\end{remark}

\noindent  Our first main result (see Theorem~\ref{th:1.1} later on)
states that  for every  $T>0$ the Cauchy problem for
system~\eqref{e1}--\eqref{bord-chi}  admits at least one solution.
In this way, we parallel the global existence result of~\cite{bbr3}:
therein, as we mentioned before, we considered a slightly different
free energy on the contact surface, which resulted in an equation
governing the evolution of the parameter $\chi$ simpler
than~\eqref{eqII}. Nonetheless, the proof of Theorem~\ref{th:1.1}
(which shall be developed in Section~\ref{rigoroso}),  closely
follows the argument developed in~\cite{bbr3}. It hinges upon a
double approximation procedure (depending on two approximating
parameters),   and a subsequent passage to the limit argument with
respect to  the mentioned parameters. One of them is used to
regularize the nonlinearities in the equations by means of Yosida
approximations. Furthermore, some viscosity terms in    $\teta$ and
$\teta_s$ are added in \eqref{e1} and \eqref{eqtetas}, depending on
the second parameter. The \emph{local} existence of a solution for
the approximate system (supplemented with suitable regularized
initial data for $\teta$ and $\teta_s$, due to the presence of
viscosity), is obtained with the Schauder theorem, while uniqueness
follows by contracting arguments. Hence,     we conclude the
existence of global-in-time solutions by proving suitable a priori
estimates (which in fact  directly hold in the time interval
$(0,+\infty)$, as they do not depend on the final time horizon $T$),
independent of the approximating parameters. The very same estimates
allow us to pass to the limit in the approximate problem,  firstly
as  the viscosity parameter and secondly as the parameter of the
Yosida regularizations
 vanish.  Finally, we point out  that, due to the strongly nonlinear character
 of system~\eqref{e1}--\eqref{bord-chi}, we do not expect
  uniqueness of solutions for  the related Cauchy problem.

\subsection*{Large-time analysis}

As previously mentioned, the ultimate aim of this paper is
investigating the large-time behaviour of  system
\eqref{e1}--\eqref{bord-chi} (supplemented with suitable initial
conditions). More precisely, we are interested in finding cluster
points of solution trajectories and characterizing a sort of
thermomechanical equilibrium of the system in the limit, in which
there is no dissipation. This corresponds to proving that solution
trajectories converge to solutions of the stationary problem
associated with our system,  in which dissipation is zero.

Now,  some results   in this direction  have been obtained in the
literature concerning the long-time behaviour of phase-field systems
with non-convex
 potentials (see, for example, \cite{colli-gilardi-laurencot-novickcohen, fs05, gps06, Krejci-Zheng,
 nuovaref}). 
 Typically, these results apply to binary systems (i.e., macroscopic
deformations are not included), see among the others \cite{BFR}
dealing with a singular entropy equation.

The main difficulties related to our analysis are due to the
singular character of the entropy equations, to the presence of
general multivalued operators on the state variables, and to the
nonlinear coupling between the equations written in the domain
$\Omega$ and the ones set in $\Gamma_c$.

\begin{remark}
\upshape
 On the other hand, the analysis of the large-time behaviour
in terms of the  global attractor for the dynamical system generated
by \eqref{e1}--\eqref{bord-chi} might also be addressed. Indeed, the
existence of the global attractor would signify that the system
dissipation is controlled in the evolution. However, for the moment
being, proving the existence of the attractor seems out of our
reach. In fact, the strongly nonlinear character of the equations
essentially prevents us from obtaining those estimates on the
solutions which would guarantee the existence of a compact and
absorbing set, for the dynamical system, in the phase space dictated
by the choice of the  initial data.
\end{remark}

Prior to addressing the large-time analysis of system
\eqref{e1}--\eqref{bord-chi}, we specify that we consider a
quadruple $(\teta,\teta_s,{\bf u},\chi)$  to be a solution
of~\eqref{e1}--\eqref{bord-chi} in $(0,+\infty)$, if
$(\teta,\teta_s,{\bf u},\chi)$ fulfils~\eqref{e1}--\eqref{bord-chi}
in the finite-time interval  $(0,T)$,  for every $T>0$. Hence, to
 perform  the asymptotic analysis on the solutions
of \eqref{e1}--\eqref{bord-chi} as time  goes to $+\infty$, we shall
rely on some further  estimates improving the solution regularity of
the existence Theorem~\ref{th:2.1}.  Only in this enhanced setting,
 shall we  obtain (see Proposition \ref{prop:2.1}) the  bounds on the
solutions (in suitable functional spaces,  on the whole half-line
$(0,+\infty)$), necessary to prove that, for every solution
trajectory,  its $\omega$-limit set  (i.e., the set of its cluster
points) is non-empty.
 These further estimates shall be
first formally derived in Section \ref{s:3}, and then made rigorous
in Section \ref{rigoroso} by performing all the related calculations
on the approximate system used for proving Theorem~\ref{th:1.1}.
That is why,
 our
asymptotic analysis   solely  applies to the solutions of
\eqref{e1}--\eqref{bord-chi} originating from the abovementioned
double approximation procedure.  Once proven that the $\omega$-limit
is non-empty,  we shall show that its elements solve the stationary
system  associated with the evolutionary problem
\eqref{e1}--\eqref{bord-chi} (see Theorem \ref{th:2.1}).   As a
by-product of this procedure, we shall see that in the limit as $t
\to \infty$ the  dissipation vanishes (cf., in particular, Remark
\ref{nodissipa}).
\paragraph{Plan of the paper.} In Section~\ref{s:2} we enlist all of
the assumptions on the problem data and state our results. A
(partially formal) proof of our Theorem~\ref{th:2.1} on the
long-time behaviour of the PDE system~\eqref{e1}--\eqref{bord-chi}
is developed in Section~\ref{s:3} and rigorously justified in
Section~\ref{rigoroso}, which also contains the proof of the global
existence Theorem~\ref{th:1.1}.

\section{Main results}
\label{s:2}
\subsection{Preliminaries}
\label{ss:2.0}
\begin{notation}
\label{not:2.1} \upshape Throughout the paper,
   given a Banach space $X$, we  shall denote by
 $_{X'}\langle\cdot,\cdot\rangle_X$ the duality pairing
between  $X'$ and $X$ itself, and by $\Vert\cdot\Vert_X$  both the
norm in  $X$ and in any power of $X$;  $\Cw([0,T];X)$ shall be the
space of the weakly continuous $X$-valued functions on $[0,T]$.

 Henceforth,  we shall suppose that
  $\Omega$ is a
bounded smooth set of $\R^3$, such that $\Gamma_c$ is a smooth
bounded domain of $\R^2$, and use the notation
$$
\begin{gathered}
H:=L^2(\Omega), \quad V:=H^1(\Omega), \quad \text{and}
\\
\bfw:= \left\{\vv \in V^3 \, : \ \vv={\bf 0} \, \hbox{ a.e. on
}\Gamma_1\right\}\,,
\end{gathered}
$$
the latter space endowed with the norm induced by $V$. We shall
denote by $\mathcal{R}$  the standard Riesz operator
\begin{equation}
\label{e:riesz}
 \mathcal{R}: V \to V' \ \ \text{given by}  \ \
 \pairing{V'}{V}{\mathcal{R}(u)}{v}:= \int_{\Omega} uv + \int_{\Omega} \nabla u \nabla
 v \quad \text{for all $u,\, v \in V$,}
\end{equation}
  and  by   $\mathcal{R}_{\Gamma_c}$ the analogously defined  Riesz operator, mapping $H^1
(\Gamma_c)$ into $(H^1 (\Gamma_c))'$.
 We shall
extensively use that
\begin{align}
& \label{continuous-embedding} V \subset L^p (\Gamma_c) \ \
\text{with a continuous (compact) embedding for} \ \ 1 \leq p \leq 4
\ \ \text{($1 \leq p <4$, resp.)},
\\
& \label{continuous-embedding-2} H^1 (\Gamma_c) \subset L^p
(\Gamma_c)\ \  \text{with a compact embedding for}\ \  1 \leq p<
\infty.
\end{align}
For notational simplicity, we shall write  $\int_{\Gamma_c}{\bf
u}{\bf v}$ ($\int_{\Gamma_2}{\bf u}{\bf v}$, respectively)  for the
duality pairing $
\pairing{(H^{-1/2}(\Gamma_c))^3}{(H^{1/2}(\Gamma_c))^3}{\bf u}{\bf
v}$ between $(H^{-1/2}(\Gamma_c))^3$ and $(H^{1/2}(\Gamma_c))^3$
(between $(H^{-1/2}(\Gamma_2))^3$ and $(H^{1/2}(\Gamma_2))^3$,
resp.). Finally, given a subset $\mathcal{O} \subset \R^N$,
$N=1,2,3,$ we shall denote by $|\mathcal{O}|$ its Lebesgue measure
and,
 for a given $v \in V'$,   the symbol $m(v)$  shall signify  its
mean value $1/|\Omega| \pairing{V'}{V}{v}{1} $.
\end{notation}
\paragraph{Variational formulation of the elasticity equation.}
We introduce the standard bilinear forms which allow to give a
variational formulation of (the boundary value problem for)
equation~\eqref{eqI}. As usual in elasticity theory, we may assume
that the material is isotropic and hence suppose that
 the rigidity matrix $K$ in
\eqref{eqI}--\eqref{condIii} can be represented as
$$
K\varepsilon({\bf u})=\lambda\hbox{tr }\left(\varepsilon({\bf
u})\right){\bf 1}+2\mu\varepsilon({\bf u}),
$$
where $\lambda,\mu>0$ are the so-called Lam\'e constants and ${\bf
1}$ is the identity matrix. Also, for the sake of simplicity but
without loss of generality, we  set $K_v={\bf 1}$ in
\eqref{eqI}--\eqref{condIii}.  Therefore, \eqref{eqI} may be
formulated by means of the following
  bilinear symmetric forms
$a, b : \bfw \times \bfw \to \R$,   defined~by
$$
\begin{aligned}
& a({\bf u},{\bf v}):=\lambda\int_{\Omega} \dive({\bf u})\dive({\bf
v})+2\mu\sum_{i,j=1}^3\int_\Omega\varepsilon_{ij}({\bf u})
\varepsilon_{ij}({\bf v})  \quad  \text{for all $\uu,\ \vv\in
\bfw$,}
\\
& b(\uu,\vv)=\sum_{i,j=1}^3\int_\Omega
\varepsilon_{ij}(\uu)\varepsilon_{ij}(\vv)
 \quad \text{for all $\uu,\ \vv\in
\bfw.$}
\end{aligned}
$$
Note that the forms $a(\cdot,\cdot)$ and $b(\cdot,\cdot)$ are
continuous and, since $\Gamma_1$ has positive measure,
 by Korn's   inequality they are $\bfw$-elliptic as well, so that
\begin{align}
& \label{form-a}  \exists\, C_a, \ K_a >0\, : \ \ a({\bf u},{\bf
u})\geq C_a\Vert{\bf u}\Vert^2_\bfw\, \qquad |a(\uu,\vv)| \leq K_a
\| \uu \|_{\bfw}\| \vv\|_{\bfw}\quad \text{for all $ \uu,\ \vv\in
\bfw$,}
\\
& \label{form-b}  \exists\, C_b, \ K_b >0\, : \ \ b({\bf u},{\bf
u})\geq C_b\Vert{\bf u}\Vert^2_\bfw\, \quad |b(\uu,\vv)| \leq K_b \|
\uu \|_{\bfw}\| \vv\|_{\bfw}\quad  \text{for all $\uu,\ \vv\in
\bfw.$}
\end{align}

\subsection{A global existence result}
\label{ss:2.1}
\paragraph{Statement of the assumptions.}

In equation~\eqref{eqII} we consider
\begin{equation}
  \label{A5}\tag{2.H1}
  \text{a maximal monotone
operator}  \ \ \beta: \R \to 2^{\R},
\end{equation}
and  denote by $\widehat\beta : \overline{\dom (\beta)} \to
(-\infty,+\infty]$ a  proper, l.s.c. and convex function  such that
$\beta=\partial\widehat\beta$. Instead, of dealing with the
pointwise operator $\partial I_-:\R\to 2^{\R}$ in~\eqref{condIii},
we shall work with a suitable generalization, defined in the duality
relation between $H^{-1/2}(\Gamma_c)^3$
 and $H^{1/2}(\Gamma_c)^3$. To this aim, we introduce
\begin{equation}
\label{hyp:alpha} \tag{2.H2}
\begin{gathered}
 \widehat \alpha : (H^{1/2}(\Gamma_c))^3  \rightarrow [0, + \infty] \, \text{a  proper,
 convex and  l.s.c. functional,}
\\
\text{with $ \widehat\alpha(\mathbf{0})=0=\min \widehat\alpha $,}
\end{gathered}
\end{equation}
and  set \[  \alpha:= \partial \widehat\alpha: \,
 (H^{1/2}(\Gamma_c))^3\to 2^{(H^{-1/2}(\Gamma_c))^3}.
\]
\begin{remark}
\upshape
 In order to render the impenetrability constraint mentioned
in the Introduction by means of the operator $\alpha$, we may
proceed as follows. We let $j({\bf u})=I_-({\bf u}\cdot{\bf n})$ and
associate with $j$ the following functional
\begin{align}\label{jh1}
&\widehat\alpha(\mathbf{v})=\int_{\Gamma_c} j({\bf
v})\quad\hbox{if}\quad {\bf v}\in
(H^{1/2}(\Gamma_c))^3\quad\hbox{and}\quad
j({\bf v})\in L^1(\Gamma_c),\\
\label{jh2}
&\widehat\alpha(\mathbf{v})=+\infty\quad\hbox{otherwise}.
\end{align}
Since  $\widehat \alpha$ is a proper, convex and lower
semicontinuous functional on  $(H^{1/2}(\Gamma_c))^3$, its
subdifferential (cf.~\cite[Cap.~II, p.~52]{barbu})
\begin{equation}\label{subvp}
\alpha:=\partial\widehat\alpha:\,(H^{1/2}(\Gamma_c))^3\to
2^{(H^{-1/2}(\Gamma_c))^3}
\end{equation}
is a maximal monotone operator. Notice that, in this case,
\eqref{jh1} implies that, if ${\bf \eta}\in\alpha({\bf v})$, then
${\bf v}$ belongs to the domain of $j$ and thus fulfils ${\bf
v}\cdot{\bf n}\leq 0$, which corresponds to the impenetrability
condition.
\end{remark}
\noindent
 We  assume that the nonlinearities $\sigma$ and $\lambda$
 comply with
\begin{equation}
\label{hyp-sig} \tag{2.H3} \sigma \in {\rm C}^{1,1} (\R)\,,
\end{equation}
\begin{equation}
\label{hyp-lambda} \tag{2.H4}
  \lambda \in {\rm C}^{1,1} (\R)\,,
\end{equation}
(and denote by $L_{\sig}$ and  $L_{\lambda}$  the Lipschitz
constants of the functions $\sig' : \R \to \R$ and $\lambda' : \R
\to \R$, respectively), and that (cf. Remark \ref{rimane})
\begin{equation}
\label{hyp-k} \tag{2.H5}
\begin{gathered}
  k \, : \R \to (0,+\infty)\ \  \text{is
Lipschitz continuous , with Lipschitz constant $L_{k}$, and}
\\
\exists\, c_k >0 \  \ \forall\, x \in \R\, : \ \ k(x) \geq c_k\,.
\end{gathered}
\end{equation}
\begin{remark}
\label{rem:consequences} \upshape We point out that
 \eqref{hyp-sig}, \eqref{hyp-lambda}, and~\eqref{hyp-k}
respectively entail that
\begin{subequations}
\begin{align}
& \label{e:crescita-c-sigma}
 \exists\, C_\sigma>0 \ \forall\, x \in \R\ :
\ \  |\sigma(x)| \leq C_\sigma (x^2 +1)\,,
\\ &
\label{e:allafineserve}
 \exists\, C_\lambda>0 \ \forall\, x \in \R\ :
\ \  |\lambda'(x)| \leq C_\lambda (|x|+1)\,,
\\
& \label{e:allafineserve-bis}
 \exists\, C_k>0 \ \forall\, x \in \R\ :
\ \  |k(x)| \leq C_k (|x|+1)\,.
\end{align}
\end{subequations}
\end{remark}
\noindent As far as  the problem data are concerned, we  suppose
\begin{equation}
\label{hypo-h} \tag{2.H6}  h \in L^2 (0,T;V')\cap L^1 (0,T;H) \,,
\end{equation}
\begin{equation}
\label{hypo-f} \tag{2.H7} \mathbf{f} \in L^2 (0,T;\mathbf{W}')\,,
\end{equation}
\begin{equation}
\label{hypo-g} \tag{2.H8} \mathbf{g} \in L^2 (0,T;
(H^{-1/2}(\Gamma_2))^3)\,.
\end{equation}
It follows from~\eqref{hypo-f}--\eqref{hypo-g} that, defining
$\mathbf{F}:(0,T) \to \bfw'$ via
$$
 \pairing{\bfw'}{\bfw}{\mathbf{F}(t)}{\vv}:=\pairing{\bf W'}{\bf W}{\mathbf{f}(t)}{\vv}
 +\int_{\Gamma_{2}} \mathbf{g}(t) \cdot {\vv}
\quad \forall\, {\vv} \in W \quad \forae \, t \in (0,T),
$$
there holds
\begin{equation}
\label{effegrande} \mathbf{F} \in L^2(0,T;\bfw') \,.
\end{equation}
Finally, we require that the initial data fulfil
\begin{align}
& \label{cond-teta-zero} \teta_0 \in L^{\bar{p}}(\Omega)\,,  \
\text{with $\bar{p} \geq \frac65$,}  \ \ \text{and} \ \ \ln(\teta_0)
\in H\,,
\\
& \label{cond-teta-esse-zero} \teta_s^0  \in
L^{\bar{q}}(\Gamma_c)\,, \ \text{with $\bar{q}>1$,} \ \ \text{and} \
\ \ln(\teta_s^0) \in L^2 (\Gamma_c)\,,
\\
& \label{cond-uu-zero} {\bf u}_0 \in \bfw \ \ \text{and} \ \  \uu_0
\in \dom (\widehat{\alpha})\,,
\\
& \label{cond-chi-zero} \chi_0 \in \hunoc  \ \ \text{and} \ \
\widehat{\beta}(\chi_0)\in
 L^1(\Gamma_c)\,.
\end{align}
Note that the first of~\eqref{cond-teta-zero} and
of~\eqref{cond-teta-esse-zero} respectively yield
\begin{equation}
\label{sto-anche-li} \teta_0 \in V', \qquad \teta_s^0 \in \hunoc'\,.
\end{equation}
\paragraph{Variational formulation and existence theorem.}
The variational formulation of the initial-boundary value problem
for system~\eqref{e1}--\eqref{bord-chi} reads as follows.
\\
{\bf Problem $\PP$.} \, Under the standing
assumptions~\eqref{A5}--\eqref{hypo-g}, given a quadruple of initial
data
 $(\teta_0, \teta_s^0 , \uu_0, \chi_0)$ complying
 with~\eqref{cond-teta-zero}--\eqref{cond-chi-zero},
find functions  $(\teta,  \teta_s,   \uu,\chi,\eeta, \xi,\zeta)$,
with the regularity
\begin{align}
& \begin{aligned} &\teta \in L^2 (0,T; V) \cap L^\infty (0, T;L^1
(\Omega))\,, \\ & \ln(\teta) \in L^\infty (0,T;H) \cap H^1
(0,T;V')\,,
\end{aligned}
 \label{reg-teta}\\
&
\begin{aligned}
&\teta_s \in  L^2 (0,T; \hunoc) \cap L^\infty (0, T;L^1
(\Gamma_c))\,, \\ &\ln(\teta_s) \in L^\infty (0,T;L^2 (\Gamma_c))
\cap H^1 (0,T;H^1(\Gamma_c)')\,,
\end{aligned}
\label{reg-teta-s}\\
&{\bf u}\in H^1(0,T;\bfw)\,,\label{reguI}\\
&  \eeta\in L^2(0,T; (H^{-1/2}(\Gamma_c))^3)\,, \label{etareg}
\\
&\chi \in L^{2}(0,T;H^2 (\Gamma_c))  \cap L^{\infty}(0,T;H^1
(\Gamma_c)) \cap H^1 (0,T;L^2 (\Gamma_c))\,, \label{regchiI}
\\
& \xi\in L^2(0,T; L^{2}(\Gamma_c))\,, \label{xireg}
\\
& \zeta \in L^\infty(\Gamma_c \times (0,T))\,, \label{zetareg}
\end{align}
 fulfilling  the initial
 conditions
\begin{align}
& \label{iniw} \teta(0)=\teta_0 \quad \aein \ \Omega\,,
\\
& \label{iniz} \teta_s(0)=\teta_s^0 \quad \aein \ \Gamma_c\,,
\\
& \label{iniu} \uu(0)=\uu_0 \quad {\aein \ \Omega}\,,
\\
& \label{inichi} \chi(0)=\chi_0 \quad {\aein \ \Gamma_c}\,,
\end{align}
  and
 \begin{align}
&
 \label{teta-weak}
\begin{aligned}
\pairing{V'}{V}{\partial_t \ln(\teta)}{v} -\int_{\Omega}
\dive(\uu_t) \, v +\int_{\Omega} &  \nabla \teta \, \nabla v +
\int_{\Gamma_c} k(\chi)  (\teta-\teta_s) v \\ & =
\pairing{V'}{V}{h}{v} \quad \forall\, v \in V \ \hbox{ a.e. in }\,
(0,T)\,,
\end{aligned}
\\
& \label{teta-s-weak}
\begin{aligned}
\pairing{{\hunoc}'}{\hunoc}{\partial_t \ln(\teta_s)}{v} &
-\int_{\Gamma_c}
\partial_t \lambda(\chi) \, v    +\int_{\Gamma_c} \nabla \teta_s  \, \nabla
v \\ &= \int_{\Gamma_c} k(\chi) (\teta-\teta_s) v  \quad \forall\, v
\in \hunoc  \ \hbox{ a.e. in }\, (0,T)\,,
\end{aligned}
\\
 &
\begin{aligned}
 b({\bf {u}}_t,\vv)+a(\uu,\vv)+
\int_{\Omega} \teta \dive (\vv)& + \int_{\Gamma_c}(\chi^+{\bf u}
+\eeta )  \cdot{\bf v} \\ &= \pairing{\bfw'}{\bfw}{\mathbf{F}}{\vv}
 \quad \forall \vv\in \bfw \ \hbox{ a.e. in }\, (0,T)\,,
 \end{aligned}
\label{eqIa}\\
&\eeta\in \alpha(\uu) \ \ \text{in $(H^{-1/2}(\Gamma_c))^3 $}\ \
\hbox { a.e. in }\, (0,T), \label{incl1}
\\
 &\chi_t -\Delta\chi+ \xi
+\sig'(\chi)=-\lambda'(\chi) \teta_s-\frac 1 2 \zeta \vert{\bf
u}\vert^2 \quad\hbox{a.e. in } \Gamma_c\times
(0,T),\label{eqIIa}\\
 &\xi
\in \beta(\chi)\ \ \hbox { a.e. in }\, \Gamma_c\times (0,T),
\label{inclvincolo}\\
 &\zeta
\in \HH(\chi)\ \ \hbox { a.e. in }\, \Gamma_c\times (0,T),
\label{inclheaviside}\\
&\partial_{\nn_s}
  \chi=0\text{ a.e. in } \partial\Gamma_c\times (0,T)\,.\label{bordo1}
  \bigskip
\end{align}
\noindent Notice that the contribution $-\lambda'(\chi)\teta_{\rm
eq}$ occurring in~\eqref{eqII}  has been incorporated into the term
$\sig'(\chi)$ in~\eqref{eqIIa}.
\begin{theorem}
\label{th:1.1}
 Assume \eqref{A5}--\eqref{hypo-g}.
\begin{enumerate}
\item Then,  Problem~$\PP$ admits a  global solution $(\teta,  \teta_s,
\uu,\chi,\eeta, \xi,\zeta)$ on the interval $(0,T)$.
\item
If,  in addition,
\begin{equation}
\label{data-further-reg}
\begin{aligned}
\mathbf{f} \in W^{1,1}(0,T; \mathbf{W}'), \qquad g \in
 W^{1,1}(0,T;(H^{-1/2}(\Gamma_2))^3), \qquad h \in W^{1,1}(0,T; V'),
\end{aligned}
\end{equation}
then there exists a solution $(\teta,  \teta_s,   \uu,\chi,\eeta,
\xi,\zeta)$ having for all $\delta>0$ the further regularity
\begin{subequations}
\label{e:further-reg}
\begin{align}
& \label{e:further-reg1} \teta \in L^\infty (\delta,T; V) \cap H^1
(\delta, T; L^{12/7}(\Omega))\,,
\\
& \label{e:further-reg2} \teta_s \in L^\infty (\delta,T; \hunoc)
\cap H^1 (\delta, T; L^{2-\rho}(\Gamma_c)) \ \ \text{for all $\rho
\in (0,2)$,}\,
\\
& \label{e:further-reg3} \chi \in L^\infty(\delta, T; H^2
(\Gamma_c)) \cap H^1 (\delta, T; \hunoc) \cap W^{1,\infty}(\delta,T;
L^2(\Gamma_c))\,,
\\
& \label{e:further-reg3-bis} \xi \in L^\infty(\delta, T; L^2
(\Gamma_c))\,,
\\
& \label{e:further-reg4} \uu \in W^{1,\infty} (\delta,T; \bfw)\,.
\end{align}
\end{subequations}
\end{enumerate}
\end{theorem}
\noindent From \eqref{e:further-reg1}--\eqref{e:further-reg3} it
follows in particular that for all $\delta>0$
 \begin{align}
 \label{e:imply}
 \teta \in \Cw ([\delta,T]; V)\,,  \qquad  \teta_s \in \Cw
([\delta,T]; H^1 (\Gamma_c))\,,\qquad \chi\in \Cw ([\delta,T];
H^2(\Gamma_c))\,.
\end{align}
We refer to Remark~\ref{rem:basta} for some further comments
concerning the above statement.
\begin{remark}
\upshape
 The regularity of $\teta_0$ and $\teta_s^0$ required in
\eqref{cond-teta-zero}--\eqref{cond-teta-esse-zero} turns out to be
necessary in the proof of our existence result Theorem~\ref{th:1.1}
for (the Cauchy problem for) system~\eqref{e1}--\eqref{bord-chi}.
Indeed, since we are going to prove existence of solutions by
passing to the limit in a viscosity approximation of equations
\eqref{e1} and \eqref{eqtetas}, we shall need to dispose of  more
regular approximate initial data. Our construction of such data (see
Lemma~\ref{dati-iniziali-lemma}) apparently hinges upon the
regularity~\eqref{cond-teta-zero}--\eqref{cond-teta-esse-zero} of
$\teta_0$ and $\teta_s^0$.

However,  we are not able to recover for $\teta$ and $\teta_s$ the
regularity corresponding to assumptions \eqref{cond-teta-zero} and
\eqref{cond-teta-esse-zero}, namely $\teta \in \mathrm{C}^0 ([0,T];
L^{\bar{p}} (\Omega))$  and $\teta_s \in \mathrm{C}^0 ([0,T];
L^{\bar{q}} (\Omega))$, with $\bar{p}$ and $\bar{q}$ as in
\eqref{cond-teta-zero}--\eqref{cond-teta-esse-zero}. This is mainly
due to the highly nonlinear character of PDE system and,  in some
sense,   to the fact that the natural initial conditions for
\eqref{teta-weak} and \eqref{teta-s-weak} are written for
$\ln(\teta)$ and $\ln(\teta_s)$. Nonetheless, notice that the
regularity required for the initial data is preserved (see
\eqref{e:imply}) for $t\geq\delta>0$, for every $\delta>0$, in the
more regular framework of \eqref{data-further-reg}.
\end{remark}
\noindent The proof of the above result is based on a double
approximation procedure which we shall detail in
Section~\ref{rigoroso}. The related passage to the limit relies on
suitable a priori estimates on the approximate solutions, which we
shall formally perform on the (un-approximated)
system~\eqref{teta-weak}--\eqref{bordo1}, and directly on the
time-interval $(0,+\infty)$, within the (formal) proof of
Proposition~\ref{prop:2.1}. Such estimates shall be rendered
rigorous in Sec.~\ref{ss:a.2}.  Finally, in Sec.~\ref{ss:a.3} we
shall conclude the proof of Theorem~\ref{th:1.1}.
\subsection{Results on the long-time behaviour of Problem $\PP$}
\label{ss:2.2} Within the scope of the present section, we shall say
that
\[
\begin{gathered}
\text{ a quadruple $(\teta,  \teta_s, \uu,\chi)$, with the
regularity~\eqref{reg-teta}--\eqref{reguI} and \eqref{regchiI},} \\
\text{ is a solution to Problem~$\PP$  if it
fulfils~\eqref{teta-weak}--\eqref{teta-s-weak} and}
\\
\text{ there exists a triple $(\eeta,\xi,\zeta)$ for
which~\eqref{etareg}, \eqref{xireg}, \eqref{zetareg}
and~\eqref{eqIa}--\eqref{bordo1} hold.}
\end{gathered}
\]

In view of the long-time analysis of the solutions to Problem~$\PP$,
we shall hereafter suppose that
\begin{equation}
\label{e:consequences} \tag{2.H9} \forall\, R > 0  \ \ \exists\,
C_{R}>0 \ \ \forall\, x \in \textrm{dom} (\widehat{\beta}) \, : \ \
\widehat{\beta}(x) + C_R \geq Rx^2\,.
\end{equation}
Notice that~\eqref{e:consequences} is trivially fulfilled in the
case $\textrm{dom} (\widehat{\beta})$ is a bounded interval,
whereas, if $\textrm{dom} (\widehat{\beta})$ is unbounded, it is
implied by a super-quadratic growth of $\widehat{\beta}$ at
infinity. We shall also require some summability on $(0,+\infty)$
for the problem data:
\begin{align}
& \label{cond-f-infty} \tag{2.H10}
 \mathbf{f} \in L^\infty (0,+\infty;\bfW')  \ \ \text{and} \ \ \mathbf{f}_t \in L^1 (0,+\infty;\bfW')\,,
\\
&
 \label{cond-g-infty}
\tag{2.H11} \mathbf{g} \in L^\infty (0,+\infty;
(H^{-1/2}(\Gamma_2))^3)\ \ \text{and} \ \  \mathbf{g}_t \in L^1
(0,+\infty; (H^{-1/2}(\Gamma_2))^3) \,,
\\
&
 \label{cond-h-infty}
 \tag{2.H12} h \in  L^\infty(0,+\infty; V') \cap L^1(0,+\infty; H)
\ \  \text{and}  \ \ h_t \in L^1(0,+\infty; V')\,.
\end{align}
The above assumptions  yield in particular
\begin{align}
& \label{cond_F-infty} \mathbf{F} \in L^\infty (0,+\infty;\bfW')  \
\ \text{and} \ \ \mathbf{F}_t \in L^1 (0,+\infty;\bfW')\,,
\\
&
 \label{e:useful-consequence} h \in L^2(0,+\infty; V')\,.
\end{align}
Furthermore,  using~\eqref{cond_F-infty}, it is not difficult to
prove (see~\cite[Remark~2.3]{bbr2} for all details) that
\begin{equation}
 \label{effe_infinito}
\exists\,\mathbf{F}_\infty\in \bfw'\,: \ \ \mathbf{F} (t) \to
\mathbf{F}_\infty \ \ \text{in $\bfw'$ as $t \to +\infty$.}
 \end{equation}

Hence, Theorem~\ref{th:1.1} ensures that for every  quadruple of
initial data $(\teta_0, \teta_s^0 , \uu_0, \chi_0)$
 fulfilling~\eqref{cond-teta-zero}--\eqref{cond-chi-zero} there
 exists (at least) a solution trajectory $(\teta,\teta_s,\uu,\chi): (0,+\infty) \to V \times \hunoc
 \times \bfW \times \hunoc $ originating from $(\teta_0, \teta_s^0 , \uu_0,
 \chi_0)$. The ensuing Proposition~\ref{prop:2.1} contains some
 suitable  large-time a priori estimates for such trajectories.
Such bounds shall enable us to conclude that the associated
 $\omega$-limit set~\eqref{omega-lim} is non-empty, and that its
 elements solve the stationary system associated with Problem~$\PP$ (see Theorem~\ref{th:2.1}).

 As we shall see, these results in fact   hold for a class of solutions of Problem
$\PP$, namely \emph{approximable solutions} which, in order not to
overburden the paper, we  shall precisely define  in
Section~\ref{rigoroso} only (cf. Definition~\ref{def:a2}). Here, we
may just mention that the notion of \emph{approximable solution} is
 tightly linked to the approximation procedure developed in
 Sec.~\ref{rigoroso}
  to prove the global existence of solutions to
Problem~$\PP$ (see Theorem \ref{th:1.1}). Such a solution notion
allows us to perform rigourously on
system~\eqref{teta-weak}--\eqref{bordo1} some of the  a priori
estimates on which our  large-time analysis relies (cf.
Remark~\ref{rem:formal}).
\paragraph{Long-time  a priori estimates.}
\begin{proposition}
\label{prop:2.1} Assume~\eqref{A5}--\eqref{hyp-k} and
\eqref{e:consequences}--\eqref{cond-h-infty}. Let $(\teta_0,
\teta_s^0, \uu_0,
 \chi_0)$ be a quadruple of initial data  complying
 with~\eqref{cond-teta-zero}--\eqref{cond-chi-zero}. Then, there
 exists a constant $K_1>0$, only depending on  the functions $\lambda$, $k$, $\sigma$,
  and on the quantity
\begin{equation}
\label{def:M}
\begin{aligned}
 M:= \|  \teta_0 \|_{L^1 (\Omega)}  + \| \teta_s^0
\|_{L^1 (\Gamma_c)} & + \| \uu_0\|_{\bfw} + \widehat{\alpha}(\uu_0)
+ \|\chi_0\|_{H^1 (\Gamma_c)} + \| \widehat{\beta}(\chi_0)\|_{L^1
(\Gamma_c)} \\  &  +
 \| \mathbf{F} \|_{L^\infty (0,+\infty;\bfw')} +
 \| \mathbf{F}_t \|_{L^1 (0,+\infty;\bfw')} \\ & + \| h
\|_{L^\infty (0,+\infty;V') \cap L^1 (0,+\infty;H)} + \| h_t \|_{L^1
(0,+\infty;V')}\,,
\end{aligned}
\end{equation}
such that for  every \emph{approximable solution} (in the sense of
Definition~\ref{def:a2})  $(\teta, \teta_s, \uu,\chi)$ to
Problem~$\PP$, fulfilling  initial
conditions~\eqref{iniw}--\eqref{inichi}, there holds
\begin{subequations}
\label{aprio-long}
\begin{align}
& \label{aprio-long1} \| \nabla \teta \|_{L^2 (0,+\infty; H)}  + \|
\nabla \teta_s \|_{L^2 (0,+\infty; L^2 (\Gamma_c))} \leq K_1\,,
\\
 &
 \label{aprio-long2}
 \| \teta \|_{L^\infty (0,+\infty; L^1 (\Omega))} + \|
\teta_s \|_{L^\infty (0,+\infty; L^1 (\Gamma_c))}\leq K_1\,,
\\
& \label{aprio-long2-bis}
 \| \teta -\teta_s \|_{L^2 (0,+\infty; L^2 (\Gamma_c))} \leq K_1\,,
\\
& \label{aprio-long3}
  \| \uu_t \|_{L^2 (0,+\infty;\bfw)} + \| \uu \|_{L^\infty
(0,+\infty;\bfw)} + \|\widehat{\alpha}(\uu)\|_{L^\infty
(0,+\infty)} +  \|\chi_t \|_{L^2 (0,+\infty; L^2 (\Gamma_c))} \leq
K_1\,,
\\
& \label{aprio-long4}
 \|
\chi \|_{L^\infty (0,+\infty;H^1 (\Gamma_c))} + \|
\widehat{\beta}(\chi)\|_{L^\infty (0,+\infty;L^1 (\Gamma_c))} \leq
K_1\,,
\\
& \label{est:2} \| \partial_t \ln(\teta) \|_{L^2 (0,+\infty; V')} +
\| \partial_t \ln(\teta_s) \|_{L^2 (0 ,+\infty; \hunoc')} \leq
K_1\,.
\end{align}
\end{subequations}
Furthermore, for all $\delta>0$
  there exist constants $K_2(\delta), \, K_{3}(\delta,\rho) >0$,
  depending on $\delta$,   on  the functions $\lambda$, $k$, $\sigma$, and on the quantity
$M$~\eqref{def:M} ($K_{3}(\delta,\rho)$ on $\rho \in (0,2)$ as
well), but independent of $(\teta, \teta_s, \uu,\chi)$, such that
 the following estimates hold
 \begin{subequations}
 \label{est-delta}
\begin{align}
& \label{est:1} \begin{aligned} \| \teta \|_{L^\infty
(\delta,+\infty;V)} &  + \| \teta_s \|_{L^\infty
(\delta,+\infty;H^1(\Gamma_c))} + \| \chi \|_{L^\infty
(\delta,+\infty;H^2 (\Gamma_c))} \\ &  + \| \chi_t\|_{L^2
(\delta,+\infty;\hunoc) \cap  L^\infty (\delta,+\infty;L^2
(\Gamma_c))} + \|\uu \|_{W^{1,\infty} (\delta,+\infty;\bfw)}   \leq
K_2(\delta)\,,
\end{aligned}
\\
 & \label{est:4}
 \|
\teta_t \|_{L^2 (\delta,+\infty;L^{12/7}(\Omega))} + \|
\partial_t\teta_s\|_{L^2 (\delta,+\infty;L^{2-\rho}(\Gamma_c))}  \leq
K_{3}(\delta,\rho) \,.
\end{align}
\end{subequations}
\end{proposition}\noindent

\begin{remark}
\label{rem:formal} \upshape In Section~\ref{ss:3.1} we shall give a
formal proof of the above result, in which all the estimates leading
to~\eqref{aprio-long} and~\eqref{est-delta} shall be performed on
the PDE system~\eqref{teta-weak}--\eqref{bordo1} directly. In
particular, this shall involve the formal differentiation of
equations~\eqref{eqIa} and~\eqref{eqIIa}, as well as  formally
testing~\eqref{teta-weak}, \eqref{teta-s-weak} by $\teta_t$ and
$\partial_t \teta_s$, respectively. All of these calculations shall
be  rigorously justified, by working on a suitable approximation of
Problem~$\PP$,  in Section~\ref{ss:a.2}.
\end{remark}
\paragraph{Results on the  $\omega$-limit of solution trajectories.}
 Now,
for  a  given quadruple of initial data $(\teta_0,\teta_s^0, \uu_0,
\chi_0)$, complying
with~\eqref{cond-teta-zero}--\eqref{cond-chi-zero}, let
 $(\teta, \teta_s, \uu,\chi)$ be an \emph{approximable solution} starting
from $(\teta_0,\teta_s^0, \uu_0, \chi_0)$, (its existence is
ensured by Theorem \ref{th:1.1}). We aim to investigate the cluster
points for large times of the trajectory of $(\teta, \teta_s,
\uu,\chi)$,
 in the topology of the space
$H^{1-\epsilon}(\Omega) \times H^{1-\epsilon}(\Gamma_c)  \times
(H^{1-\epsilon}(\Omega))^3 \times  H^{2-\epsilon}(\Gamma_c) $,  with
an arbitrary $\epsilon >0$.
 To this aim, we
define the $\omega$-limit set $\omega(\teta,\teta_s,\uu, \chi)$ of
the trajectory $(\teta(t),\teta_s(t),\uu(t), \chi(t))_{t \geq 0}$ as
\begin{equation}
\label{omega-lim}
\begin{aligned}
\omega(\teta,\teta_s,\uu, \chi):= & \big \{
(\teta_\infty,\tetasinf,\uu_\infty, \chi_\infty) \in V \times \hunoc
\times  \bfw \times H^1(\Gamma_c)
 \, : \\ &  \exists
\{t_n \} \subset [0,+\infty),  \ t_n \nearrow +\infty \  \text{as $n
\uparrow \infty$},
\\
&\text{with} \ (\teta(t_n),\teta_s(t_n), \uu(t_n), \chi(t_n)) \to
(\teta_\infty,\tetasinf,\uu_\infty, \chi_\infty) \\ &  \text{in $
H^{1-\epsilon}(\Omega) \times H^{1-\epsilon}(\Gamma_c)  \times
\left(H^{1-\epsilon}(\Omega)\right)^3 \times
H^{2-\epsilon}(\Gamma_c) $} \big\}.
\end{aligned}
\end{equation}
For simplicity, we choose to omit the dependence  on the initial
data $(\teta_0,\teta_s^0, \uu_0, \chi_0)$  (and on the parameter
$\epsilon$), in the notation $\omega(\teta,\teta_s,\uu, \chi)$.
 Notice that the latter would be replaced by the more customary
$\omega(\teta_0,\teta_s^0, \uu_0, \chi_0)$ if we additionally
disposed of a uniqueness result for the \emph{approximable
solutions} to Problem~$\PP$.

The following theorem shall be proved in Section~\ref{ss:3.2}.
\begin{theorem}
\label{th:2.1}
 Assume~\eqref{A5}--\eqref{hyp-k} and
 \eqref{e:consequences}--\eqref{cond-h-infty}.

 Then, for every quadruple of initial data
$(\teta_0,\teta_s^0,\uu_0,\chi_0)$ complying
with~\eqref{cond-teta-zero}--\eqref{cond-chi-zero} and for every
approximable solution $(\teta(t),\teta_s(t),\uu(t), \chi(t))_{t \geq
0}$ originating from $(\teta_0,\teta_s^0,\uu_0,\chi_0)$, for all
$\epsilon>0$ the associated $\omega$-limit set
$\omega(\teta,\teta_s,\uu, \chi)$ is a non-empty, compact and
connected subset of $ H^{1-\epsilon}(\Omega) \times
H^{1-\epsilon}(\Gamma_c) \times
\left(H^{1-\epsilon}(\Omega)\right)^3 \times
H^{2-\epsilon}(\Gamma_c) $.

Furthermore,   every quadruple $(\teta_\infty,\tetasinf,\uu_\infty,
\chi_\infty) \in \omega(\teta,\teta_s,\uu, \chi)$ is a solution of
the stationary system associated with~Problem~$\PP$, namely
\begin{subequations}
 \label{e:stationary-system}
\begin{align}
& \label{cuore} \exists\,\overline{\teta}_\infty \geq 0\, : \quad
\teta_\infty(x) \equiv \overline{\teta}_\infty  \ \ \text{for a.e.
$x \in \Omega$,}\quad \tetasinf(x) \equiv \overline{\teta}_\infty  \
\  \text{for a.e. $x \in \Gamma_c$,}
\\
& \label{e:stat1}
\begin{aligned}
& a(\uu_\infty,\vv)+ \overline{\teta}_\infty \int_{\Omega} \dive
(\vv) + \int_{\Gamma_c}(\chi_\infty^+\uu_\infty +\eeta_\infty )
\cdot{\bf v}=\mathbf{F}_\infty
 \ \  \forall \vv\in \bfw\,,
 \\
  &
  \eeta_\infty \in
 \alpha(\uu_\infty) \ \ \text{in $(H^{-1/2}(\Gamma_c))^3$,}
 \end{aligned}
 \\
 &
 \label{e:stat2}
 \chi_{\infty} \in H^2 (\Gamma_c) \ \text{and} \
\begin{cases}
& -\Delta \chi_{\infty} + \xi_{\infty}+\sigma'(\chi_\infty)
=-\lambda'(\chi_\infty)\overline{\teta}_\infty  -\frac12
\zeta_\infty|\uu_\infty|^2 \quad \aein \ \Gamma_c,
\\
& \text{$\xi_\infty \in L^2 (\Gamma_c)$, \quad $ \xi_\infty \in
\beta(\chi_\infty)  \ \aein \ \Gamma_c,$}
\\
& \zeta_\infty \in L^\infty (\Gamma_c), \quad \zeta_\infty \in
H(\chi_\infty)\ \aein \ \Gamma_c,
\\
&
\partial_{\nn_s} \chi_\infty = 0 \quad \aein \  \partial\Gamma_c\,.
\end{cases}
\end{align}
\end{subequations}
\end{theorem}
\begin{remark}\label{nodissipa}
\upshape Let us emphasize that, in  the limit as $t \to \infty$, the
system is in a state of thermomechanical equilibrium. Indeed,   the
dissipation, described by the pseudo-potentials \eqref{pseudoO} and
\eqref{pseudoG}, has vanished in \eqref{cuore}--\eqref{e:stat2}.
\end{remark}
\noindent  No uniqueness result is available on the stationary
system \eqref{cuore}--\eqref{e:stat2}. Hence, one cannot deduce
directly from Theorem \ref{th:2.1} that $\omega(\teta,\teta_s,\uu,
\chi)$ is a singleton and that the \emph{whole}  solution trajectory
$(\teta,\teta_s,\uu, \chi)$ thus converges, as $t \to +\infty$, to a
unique equilibrium. However, the next result (whose proof is
postponed to  Section~\ref{ss:3.2}) shows that, under more specific
assumptions on the operator $\beta$ and on the nonlinearities
$\lambda$ and $\sigma$, it is possible to uniquely determine the
$\chi$-component of the elements in $\omega(\teta,\teta_s,\uu,
\chi)$.
\begin{corollary}
\label{coro:particular-case} Under assumptions
\eqref{A5}--\eqref{hyp-k} and
\eqref{cond-f-infty}--\eqref{cond-h-infty}, suppose further that
 \begin{subequations}
 \label{e:more-specific}
 \begin{align}
& \label{e:spec1} \beta= \partial I_{[m_*,m^*]}, \quad \text{for
some $-\infty<m_*<m^*<+\infty$,}
\\
& \label{e:spec2} \lambda  \quad \text{is non-decreasing on
$[m_*,m^*]$,}
\\
& \label{e:spec3} \sigma'(x) >0 \quad \text{for all $x \in
[m_*,m^*]$.}
 \end{align}
 \end{subequations}
Then, for any quadruple $(\teta_\infty,\tetasinf,\uu_\infty,
\chi_\infty) \in \omega(\teta,\teta_s,\uu, \chi)$ there holds
\begin{equation}
\label{chii-infi-det} \chi_\infty (x) \equiv m_* \quad  \forall \, x
\in \Gamma_c\,,
\end{equation}
and we have as $t \to +\infty$
\[
\chi(t) \to \chi_\infty  \quad \text{in $H^{2-\epsilon}(\Gamma_c) $
for all $\epsilon>0$.}
\]
\end{corollary}
\begin{remark}
\upshape Corollary \ref{coro:particular-case} ensures  that, in the
case when the latent heat is positive, $\beta$ has a bounded domain
(which is the interesting case from a physical point of view), and
cohesion in the material (which is included in the decreasing part
of $\sigma$, see \eqref{freeG}) is not too large with respect to the
remaining part of the potential $\sigma$, the glue tends to be
completely damaged in  the  large-time limit. We remark that  this
is the result one would expect from experience.
\end{remark}
\section{Proofs}
\label{s:3}
\begin{notation} \upshape
 Henceforth,  for the sake of notational
simplicity, we shall write $\langle \cdot, \cdot \rangle$  for all
the duality pairings
 $\pairing{\bfw'}{\bfw}{\cdot}{\cdot}$,
$\pairing{V'}{V}{\cdot}{\cdot}$, and
$\pairing{{\hunoc}'}{\hunoc}{\cdot}{\cdot}$,  and, further,
 denote
by the symbols
\[
\text{$c$, $c'$ $C$, $C'$ most of the (positive) constants occurring
in calculations and estimates.}
\]
Notice that the above constants shall not depend on the final
time $T$.
\end{notation}
\subsection{Proof of Proposition~\ref{prop:2.1}}
\label{ss:3.1} \noindent The following is a  \emph{formal} proof,
based on two
 main a priori estimates on system~\eqref{teta-weak}--\eqref{bordo1},
 which shall be revisited in Section~\ref{ss:a.2}.
 Therein, we shall also rigorously  prove the further solution
 regularity~\eqref{e:further-reg1}--\eqref{e:further-reg4}.
\paragraph{First (formal) estimate.} We test~\eqref{teta-weak} by
$\teta$, \eqref{teta-s-weak} by $\teta_s$, \eqref{eqIa} by $\uu_t$
and~\eqref{eqIIa} by $\chi_t$, add the resulting relations and
integrate them on the interval $(0,t)$, with $t \in (0,+\infty)$.
Now, we take into  account the \emph{formal} identities
\begin{equation}
\label{e:formal1}
\begin{aligned}
& \int_0^t \langle \partial_t \ln(\teta), \teta \rangle = \|
\teta(t) \|_{L^1 (\Omega)} - \| \teta_0 \|_{L^1 (\Omega)}\,,
\\
& \int_0^t \langle \partial_t \ln(\teta_s), \teta_s \rangle = \|
\teta_s(t) \|_{L^1 (\Gamma_c)} - \| \teta_s^0 \|_{L^1 (\Gamma_c)}\,,
\end{aligned}
\end{equation}
and the chain rule  for the  convex functionals $\widehat{\alpha}$,
 and $\widehat{\beta}$ (cf.
with~\cite[Lemma~4.1]{colli92} and~\cite[Lemma~3.3]{brezis73},
respectively), and for the smooth function $\sigma$, yielding
\begin{equation}
\label{e:rig1}
\begin{gathered}
 \int_0^t \int_{\Gamma_c} \eeta \cdot \uu_t =
\widehat{\alpha}(\uu(t)) - \widehat{\alpha} (\uu_0)\,, \qquad
\int_0^t \int_{\Gamma_c} \xi  \chi_t =
\int_{\Gamma_c}\widehat{\beta}(\chi(t)) -
\int_{\Gamma_c}\widehat{\beta} (\chi_0)\,,
\\
 \int_0^t
\int_{\Gamma_c} \sigma'(\chi)\chi_t =
\int_{\Gamma_c}\sigma(\chi(t))- \int_{\Gamma_c}\sigma(\chi_0)\,.
\end{gathered}
\end{equation}
In the same way, an integration by parts and the chain rule for
$\mathrm{p}(\cdot)=(\cdot)^+$ lead to
\begin{equation}
\label{e:to-be-cited}
\begin{aligned}
 \int_0^t
\int_{\Gamma_c} \chi^+ \uu \uu_t &  = \frac12\int_0^t
\int_{\Gamma_c} \chi^+
\partial_t |\uu|^2 \\ & = -\frac12 \int_0^t \int_{\Gamma_c}\zeta \chi_t
|\uu|^2+ \frac12\int_{\Gamma_c} \chi^+(t) |\uu(t)|^2 -
\frac12\int_{\Gamma_c} \chi_0^+|\uu_0|^2\,,
\end{aligned}
\end{equation}
where we recall that $\zeta\in \HH(\chi)=\partial \mathrm{p}(\chi)$.
 Finally, we observe
that,   by the properties~\eqref{form-a} and~\eqref{form-b}  of  the
forms $a$ and $b$, respectively, we have
\begin{equation}
\label{e:coercive-forms} \int_0^t a(\uu,\uu_t) \geq \frac{C_a}2 \|
\uu(t) \|_{\bfw}^2 - \frac12 K_a \|\uu_0\|_{\bfw}^2\,, \qquad
 \int_0^t b
(\uu_t,\uu_t) \geq C_b \int_0^t \| \uu_t \|_{\bfw}^2\,.
\end{equation}
Collecting~\eqref{e:formal1}--\eqref{e:coercive-forms} and observing
that some terms cancel out,
 we get
\begin{equation}
\label{e:3.5}
\begin{aligned}
 \| \teta(t) \|_{L^1 (\Omega)}  & + \int_0^t \| \nabla \teta \|_H^2 +
\int_0^t \int_{\Gamma_c} k(\chi) (\teta-\teta_s)^2 + \| \teta_s(t)
\|_{L^1 (\Gamma_c)} + \int_0^t \| \nabla \teta_s \|_{L^2
(\Gamma_c)}^2 \\ & + C_b  \int_0^t \| \uu_t \|_{\bfw}^2  +
\frac{C_a}2 \| \uu(t) \|_{\bfw}^2 + \widehat{\alpha} (\uu(t)) +
\frac12\int_{\Gamma_c} \chi^+(t) |\uu(t)|^2 \\ & + \int_0^t \|
\chi_t \|_{L^2 (\Gamma_c)}^2 + \frac12 \| \nabla \chi(t) \|_{L^2
(\Gamma_c)}^2 + \int_{\Gamma_c}\left( \widehat{\beta}(\chi(t)) +
\sigma(\chi(t)) \right) \\ &
\begin{aligned}
 \leq \| \teta_0\|_{L^1 (\Omega)}  & + \|
\teta_s^0 \|_{L^1 (\Gamma_c)}   + \frac{K_a}2 \| \uu_0 \|_{\bfw}^2 +
\widehat{\alpha} (\uu_0) +\frac12 \| \nabla \chi_0 \|_{L^2
(\Gamma_c)}^2\\ + & c \| \chi_0\|_{L^2 (\Gamma_c)} \| \uu_0
\|_{\bfw}^2 + \| \widehat{\beta}(\chi_0) \|_{L^1 (\Gamma_c)}+ \|
\sigma(\chi_0) \|_{L^1 (\Gamma_c)} + I_1 + I_2\,,
\end{aligned}
\end{aligned}
\end{equation}
where have also used the continuous
embedding~\eqref{continuous-embedding} for the third term on the
right-hand side of~\eqref{e:to-be-cited}. Then, we estimate
\begin{align}
 \label{e2.1.1}
\begin{aligned}
 I_1 = \int_0^t \langle h, \teta \rangle  &  \leq
\int_0^t \| h \|_{V'} \| \teta -m(\teta) \|_V + \int_0^t \| h
\|_{V'} \| m(\teta) \|_V \\ & \leq C \int_0^t  \| h \|_{V'}^2 +
 \frac12 \int_0^t \| \nabla \teta
\|_H^2 + \frac1{|\Omega|^{1/2}} \int_0^t \| h \|_{V'} \|
\teta\|_{L^1 ((\Omega)}\,,
\end{aligned}
\end{align}
the second inequality due to the Poincar\'e inequality for functions
with zero mean value, and, with an integration  by parts, we get
\begin{align}
 \label{e2.1.2}
\begin{aligned}
  I_2 = \int_0^t \langle \mathbf{F}, \uu_t \rangle
& = - \int_0^t \langle \mathbf{F}_t, \uu \rangle + \langle
\mathbf{F}(t), \uu(t) \rangle -\langle \mathbf{F}(0), \uu_0 \rangle
\\ & \leq \int_0^t \|\mathbf{F}_t\|_{\bfw'} \|\uu \|_{\bfw} +
\frac{C_a}{4}\|\uu(t) \|_{\bfw}^2 + \frac12\| \uu_0 \|_{\bfw}^2 +
C\| \mathbf{F}\|_{L^\infty (0,+\infty;\bfw')}^2\,.
\end{aligned}
\end{align}

Hence,
 taking
into account our assumptions on the
data~\eqref{cond-f-infty}--\eqref{cond-h-infty} (which
yield~\eqref{cond_F-infty} and~\eqref{e:useful-consequence}), as
well as~\eqref{cond-teta-zero}--\eqref{cond-chi-zero},   we may
apply a variant of  the Gronwall Lemma
(cf.~\cite[Lemma~A.5]{brezis73}) and conclude
estimates~\eqref{aprio-long1}, \eqref{aprio-long2},
and~\eqref{aprio-long3}, while~\eqref{aprio-long2-bis} follows from
the bound for $(k(\chi))^{1/2} (\teta -\teta_s)$ in $L^2
(0,+\infty;L^2 (\Gamma_c))$ and the fact that $k$ is bounded from
below by a strictly positive constant, see~\eqref{hyp-k}. Finally,
thanks to~\eqref{e:consequences} and~\eqref{e:crescita-c-sigma}, we
deduce from~\eqref{e:3.5} that
\begin{equation}
\label{beta-sigma} \exists\, C, \,c, \,  c'>0 \ \ \forall\, t\in
(0,+\infty)\, : \ \  c \| \chi(t)\|_{L^2 (\Gamma_c)}^2 - c' \leq
\int_{\Gamma_c}\left( \widehat{\beta}(\chi(t)) + \sigma(\chi(t))
\right)  \leq C\,.
\end{equation} Joint with the bound for $\nabla \chi$ in $L^\infty
(0,+\infty; L^2(\Gamma_c))$,  this yields  the estimate for
$\|\chi\|_{L^\infty (0,+\infty;\hunoc)}$. A fortiori,  in view
of~\eqref{e:consequences} we also recover the bound for
$\widehat{\beta}(\chi)$, and
 \eqref{aprio-long4} ensues. In the end, estimate
\eqref{est:2}  for $\partial_t \ln(\teta)$ and for $\partial_t
\ln(\teta_s)$ follows from a comparison in
equations~\eqref{teta-weak} and~\eqref{teta-s-weak}, respectively.
For example, using~\eqref{hyp-k} and~\eqref{aprio-long2-bis},
\eqref{aprio-long4} one observes that
\[
\| k(\chi) (\teta-\teta_s) \|_{L^2 (0, +\infty; V')} \leq C (\|
\chi\|_{L^\infty (0,+\infty; L^4 (\Gamma_c))} + 1) \| \teta-\teta_s
\|_{L^2 (0,+\infty; L^2(\Gamma_c))} \leq C'\,.
\]
Hence, in view of the bound~\eqref{aprio-long1} for $\nabla \teta$
in $L^2 (0,+\infty; H)$, of~\eqref{aprio-long3} for $\dive(\uu_t)$
in $L^2 (0,+\infty; H^3)$, and of~\eqref{e:useful-consequence}, one
concludes the estimate for $\partial_t  \ln(\teta) $ in $L^2
(0,+\infty; V')$.
\paragraph{Second (formal)  estimate.}
In what follows, we shall formally treat the maximal monotone
operators $\alpha$, $\HH$, and $\beta$ as nondecreasing and
\emph{Lipschitz} functions. Indeed, the following estimates can be rigorously justified as in Section
\ref{rigoroso}
regularizing these nonlinearities by their Yosida approximations (which are in fact Lipschitz functions). Further,
we shall work with a solution
$(\teta,\teta_s,\uu,\chi)$ enjoying the further
regularity~\eqref{e:further-reg} (see also \eqref{e:imply}).
Formally proceeding, for all $\delta>0$ we let a
constant $K_\delta$, only depending on $\lambda$, $k$, $\sigma$, on
the quantity $M$~\eqref{def:M}, and possibly on $\delta$, such that
\begin{equation}
\label{e:brutta} \| \teta(\delta) \|_{V}^2 + \| \teta_s
 (\delta)\|_{H^1(\Gamma_c)}^2 + \| \uu_t (\delta) \|_{\bfw}^2 + \| \chi_t
(\delta) \|_{L^2 (\Gamma_c)}^2 \leq K_\delta\,.
\end{equation}
This formal assumption shall be  discarded once we put forth the
rigorous arguments of Section \ref{rigoroso}, see
Remark~\ref{rem:ripensamento} for further comments. Hence, we are in
the position of performing the following
 calculations.
  We test~\eqref{teta-weak} by $ \teta_t$ and
\eqref{teta-s-weak} by $\partial_t \teta_s$,
differentiate~\eqref{eqIa} w.r.t. time and test it by $\uu_t$, and
differentiate~\eqref{eqIIa} w.r.t. time  and multiply it by
$\chi_t$.
 We add the resulting relations and integrate them
on the interval $(\delta,t)$, with $t \in (\delta,+\infty)$, also
adding $1/2 (\|\teta(t) \|_{L^1(\Omega)}^2 +\| \teta_s(t) \|_{L^1
(\Gamma_c)}^2 )$ to both sides. Indeed, the Poincar\'e inequality
for the zero mean value functions yields
\begin{equation}
\label{e:poincare}
\begin{aligned}
& C_P \| \teta(t) \|_V^2 \leq  \frac12 \left(\| \nabla \teta (t)
\|_{H}^2+ \| \teta(t) \|_{L^1 (\Omega)}^2\right) \\ & C_P \|
\teta_s(t) \|_{\hunoc}^2 \leq  \frac12 \left(\| \nabla \teta_s (t)
\|_{L^2(\Gamma_c)}^2+ \| \teta_s(t) \|_{L^1 (\Gamma_c)}^2\right)
\end{aligned}
\end{equation}
for some positive constant $C_P$ independent of $t \in (0,+\infty)$.
We also notice that, for some other constant $c$ also depending on
the
embeddings~\eqref{continuous-embedding}--\eqref{continuous-embedding-2},
 there holds
\begin{equation}
\label{e:l4} c  \| \teta(t) - \teta_s(t) \|_{L^4 (\Gamma_c)}^2 \leq
\frac{C_P}2 \left(\| \teta(t) \|_V^2 + \| \teta_s(t) \|_{\hunoc}^2
\right) \quad \text{for all $t \in (0,+\infty)$.}
\end{equation}
Further, we remark  that
\begin{equation}
\label{e:lampada}
\begin{aligned}
\int_\delta^t \int_{\Gamma_c} \partial_t (\chi^+ \uu) \uu_t =
\int_\delta^t \int_{\Gamma_c} \chi^+ |\uu_t|^2 +  \int_\delta^t
\int_{\Gamma_c} \HH(\chi)\chi_t \uu \uu_t\,.
\end{aligned}
\end{equation}
 Taking into account the
cancellation of some terms and
 the coercivity  and continuity of the forms
$a$ and $b$~\eqref{form-a}--\eqref{form-b},   using the
\emph{formal} identities
\begin{equation}
\label{e:formal2}
\begin{aligned}
 \int_\delta^t \langle \partial_t \ln(\teta),  \teta_t \rangle =
\int_\delta^t \int_{\Omega}\frac{|\teta_t|^2}{\teta}\,, \qquad
\int_\delta^t \langle \partial_t \ln(\teta_s),
\partial_t\teta_s \rangle =  \int_\delta^t \int_{\Gamma_c}
\frac{|\partial_t\teta_s|^2}{\teta_s}\,,
\end{aligned}
\end{equation}
as well as~\eqref{e:poincare}--\eqref{e:lampada},  we end up with
\begin{equation}
\label{e:3.6}
\begin{aligned}
& \int_\delta^t \int_{\Omega} \frac{|\teta_t|^2}{\teta} +
\frac{C_P}2 \| \teta(t) \|_V^2
 + \int_{\delta}^t \int_{\Gamma_c} k(\chi)
(\teta-\teta_s)
\partial_t(\teta-\teta_s) + \int_\delta^t \int_{\Gamma_c}
\frac{|\partial_t\teta_s|^2}{\teta_s}  \\ & +\frac{C_P}2 \|  \teta_s
(t) \|_{\hunoc}^2  + c  \| \teta(t) - \teta_s(t) \|_{L^4
(\Gamma_c)}^2 +\frac{C_b}2 \| \uu_t (t) \|_{\bfw}^2 +C_a
\int_\delta^t \|
\uu_t \|_{\bfw}^2  \\
&
\begin{aligned}
+ \int_\delta^t \int_{\Gamma_c} \chi^+ |\uu_t|^2  & + \int_\delta^t
\langle \alpha'(\uu) \uu_t, \uu_t \rangle   +\frac12 \|\chi_t (t)
\|_{L^2 (\Gamma_c)}^2  + \int_\delta^t \int_{\Gamma_c} |\nabla
\chi_t|^2\\
&   + \int_\delta^t \int_{\Gamma_c} \beta'(\chi) |\chi_t|^2 \leq C
K_\delta+ K_1^2 + I_3 + I_4 + I_5 + I_6 + I_7 + I_8 \,,
\end{aligned}
\end{aligned}
\end{equation}
in which we have controlled the term $1/2 (\|\teta(t)
\|_{L^1(\Omega)}^2 +\| \teta_s(t) \|_{L^1 (\Gamma_c)}^2 )$ on the
right-hand side by~\eqref{aprio-long2}.
 Integrating by parts, we have
\begin{align}
 \label{3.6.1}
\begin{aligned}
 I_3 = \int_\delta^t \langle h, \teta_t \rangle = &
- \int_\delta^t \langle h_t, \teta \rangle  + \langle
h(t),\teta(t)\rangle  - \langle h(\delta),\teta(\delta)\rangle\\ &
\leq \int_\delta^t \| h_t \|_{V'} \|\teta\|_V + C \| h \|_{L^\infty
(0,+\infty;V')}^2   + \frac{C_P}4 \|\teta(t) \|_{V}^2 + K_\delta\,,
\end{aligned}
\end{align}
while  we estimate
 \begin{align}
 &
  \label{3.6.2}
  I_4 =
\int_\delta^t \langle \mathbf{F}_t, \uu_t \rangle \leq \int_\delta^t
\| \mathbf{F}_t\|_{\bfw'} \| \uu_t \|_{\bfw} \,,
\\
 \label{3.6.21}
&
\begin{aligned}
 I_5= -2 \int_\delta^t \int_{\Gamma_c} &  \HH(\chi)\chi_t \uu \uu_t   \leq
2\int_\delta^t \|\chi_t\|_{L^2 (\Gamma_c)} \|\uu\|_{L^4 (\Gamma_c)}
\|\uu_t\|_{L^4 (\Gamma_c)} \\ & \leq C \| \uu \|_{L^\infty
(0,+\infty;\bfW)}  \left(\int_\delta^t \int_{\Gamma_c}
\|\chi_t\|_{L^2 (\Gamma_c)}^2 + \int_\delta^t \int_{\Gamma_c}
\|\uu_t\|_{\bfW}^2 \right) \leq C\,,
\end{aligned}
\\
 \label{3.6.22}
 &
I_6 = -\frac12\int_\delta^t \int_{\Gamma_c} \HH'(\chi)|\chi_t|^2
|\uu|^2 \leq 0\,,
\end{align}
\eqref{3.6.21} due to the continuous
embedding~\eqref{continuous-embedding} and to
estimate~\eqref{aprio-long3}, and \eqref{3.6.22} following by
monotonicity. Moreover, recalling~\eqref{hyp-sig},
\eqref{hyp-lambda}, and~\eqref{aprio-long3}, and applying the
H\"older inequality
\begin{align}
&
 \label{3.6.2.bis}
  I_7 =-
\int_\delta^t \int_{\Gamma_c} \sigma{''}(\chi) |\chi_t|^2 \leq
L_\sigma K_1^2\,,
\\
& \label{3.6.3}
 \begin{aligned}
  I_8 =-
\int_\delta^t \int_{\Gamma_c} \lambda{''}(\chi) |\chi_t|^2 \teta_s &
\leq L_\lambda \int_\delta^t \|\chi_t \|_{L^2 (\Gamma_c)} \|\chi_t
\|_{L^4 (\Gamma_c)} \|\teta_s \|_{L^4 (\Gamma_c)}\\ &  \leq \frac14
\int_\delta^t \| \nabla \chi_t\|_{L^2 (\Gamma_c)}^2 +\frac{K_1^2}4 +
C \int_\delta^t \|\chi_t \|_{L^2 (\Gamma_c)}^2 \|\teta_s
\|_{\hunoc}^2 \,,
\end{aligned}
\end{align}
$C$
 also depending on the constant of the continuous
 embedding~\eqref{continuous-embedding-2}.
 Finally, we integrate by parts the third term on the left-hand
side of~\eqref{e:3.6}, so that
\begin{equation}
\label{e:inte-parts}
\begin{aligned}
\int_{\delta}^t  & \int_{\Gamma_c} k(\chi) (\teta-\teta_s)
\partial_t(\teta-\teta_s)\\ &  =
\frac12\int_{\Gamma_c} k(\chi(t)) |\teta(t) -\teta_s(t)|^2
-\frac12\int_{\Gamma_c} k(\chi(\delta)) |\teta(\delta)
-\teta_s(\delta)|^2- I_9 \\
& \geq \frac{c_k}2 \| \teta(t) -\teta_s(t)\|_{L^2 (\Gamma_c)}^2 - c
\left( \|\chi(\delta)\|_{L^2 (\Gamma_c)} +1\right)\left(
\|\teta(\delta)\|_{V}^2 + \|\teta_s (\delta)\|_{\hunoc}^2\right) -
I_9
\\
&  \geq \frac{c_k}2 \| \teta(t) -\teta_s(t)\|_{L^2 (\Gamma_c)}^2 -
C\left( {K_\delta} +1\right)   - I_9 \,,
\end{aligned}
\end{equation}
the latter inequality due~\eqref{hyp-k},
to~\eqref{e:allafineserve-bis}, \eqref{aprio-long4}, and
\eqref{e:brutta}, while, using~\eqref{hyp-k}
and~\eqref{aprio-long3}, we estimate
\begin{equation}
\label{3.6.4}
\begin{aligned}
 I_9 &  =  \frac12 \int_{\delta}^t \int_{\Gamma_c}
k'(\chi)\chi_t|\teta-\teta_s|^2 \\ & \leq \frac{L_k}2
\int_{\delta}^t \|\chi_t \|_{L^4 (\Gamma_c)} \| \teta-\teta_s\|_{L^2
(\Gamma_c)} \|\teta-\teta_s \|_{L^4 (\Gamma_c)} \\ & \leq \frac12
\int_\delta^t \| \nabla \chi_t\|_{L^2 (\Gamma_c)}^2 +\frac{K_1^2}2 +
C \int_\delta^t \| \teta-\teta_s\|_{L^2 (\Gamma_c)}^2
\|\teta-\teta_s \|_{L^4 (\Gamma_c)}^2
\end{aligned}
\end{equation}
Next, we collect~\eqref{e:3.6}--\eqref{3.6.4}, observing that the
ninth term on the left-hand side of~\eqref{e:3.6} is non-negative,
 and so
are the tenth and the thirteenth terms, by monotonicity of $\alpha$
and $\beta$, respectively. Then, taking into account the summability
properties~\eqref{cond_F-infty} and~\eqref{e:useful-consequence}, as
well as estimate~\eqref{aprio-long3} for  $\|\chi_t
\|_{L^2(0,+\infty; L^2 (\Gamma_c))}$ and~\eqref{aprio-long2-bis} for
$\| \teta-\teta_s \|_{L^2(0,+\infty; L^2 (\Gamma_c))}$,
 we
apply the standard Gronwall lemma and  its abovementioned variant
(cf.~\cite[Lemmas A.3, A.5]{brezis73}) and conclude
estimate~\eqref{est:1} for $\teta$, $\teta_s$, $\uu$, and $\chi_t$.
A comparison in~\eqref{eqIIa} yields (possibly for a larger
$K_2(\delta)$)
\[
\| -\Delta \chi + \xi\|_{L^\infty(\delta, +\infty; L^2
(\Gamma_c))}^2 \leq K_2(\delta)\,,
\]
whence an estimate both for $\| \xi\|_{L^\infty(\delta, +\infty; L^2
(\Gamma_c))}$ and for $  \|-\Delta \chi  \|_{L^\infty(\delta,
+\infty; L^2 (\Gamma_c))}$ by the  monotonicity of $\beta$. Then,
the bound for $ \| \chi  \|_{L^\infty(\delta, +\infty; H^2
(\Gamma_c))}$ follows from elliptic regularity.
 Furthermore, the estimate on the first and
on the fourth integral terms on the left-hand side of~\eqref{e:3.6}
gives (notice that $\teta>0$ a.e. in $\Omega \times (0,+\infty)$ and
$\teta_s>0$ a.e. in $\Gamma_c \times (0,+\infty)$)
\begin{equation}
\label{e:further-bound} \| \partial_t \teta^{1/2} \|_{L^2
(\delta,+\infty; H)} +\| \partial_t \teta_s^{1/2} \|_{L^2
(\delta,+\infty; L^2(\Gamma_c))} \leq C(\delta)\,.
 \end{equation}
On the other hand, due to the continuous embeddings $V \subset L^p
(\Omega)$ for all $ 1 \leq p \leq 6$ and $\hunoc \subset L^q
(\Gamma_c)$ for all $1 \leq q <\infty$, \eqref{est:1} in particular
yields
\begin{equation}
\label{e:further-bound-1} \forall\,1 \leq q <\infty \ \ \exists\,
C_q(\delta)>0 \, : \  \ \|\teta^{1/2}\|_{L^\infty (\delta, +\infty;
L^{12}(\Omega))}+ \|\teta_s^{1/2}\|_{L^\infty (\delta, +\infty;
L^{q}(\Gamma_c))} \leq C_q(\delta)\,.
\end{equation}  Hence, \eqref{est:4} follows from~\eqref{est:1},
\eqref{e:further-bound}, \eqref{e:further-bound-1}, and the H\"older
inequality,  giving
\begin{equation}
\label{e:holder}
\begin{aligned}
&  \| \teta_t \|_{L^{12/7}(\Omega)} \leq \|
\partial_t \teta^{1/2} \|_{H} \| \teta^{1/2}\|_{L^{12}(\Omega)}
\\
&
 \| \partial_t \teta_s \|_{L^{2-\rho}(\Gamma_c)} \leq \|
\partial_t \teta_s^{1/2} \|_{L^2(\Gamma_c)} \|
\teta_s^{1/2}\|_{L^{(4-2\rho)/\rho}(\Gamma_c)} \ \  \ \text{for all
$0 < \rho < 2$.}
\end{aligned}
\end{equation}
This concludes the proof. \fin
\begin{remark}
\label{rem:ripensamento} \upshape Following~\cite{bcfg2},  we point
out that a possible way to perform estimate~\eqref{e:3.6} more
rigorously, without
 assuming~\eqref{e:brutta}, would be to
fix a smooth ``cut-off'' function $\varsigma: [0,+\infty) \to
[0,+\infty)$,  with $\varsigma, \, \varsigma' \in L^\infty
(0,+\infty)$ (for example, $\varsigma(t) = \tanh(t)$ for all $t \geq
0$), and
 test~\eqref{teta-weak} by $\varsigma \teta_t$,
\eqref{teta-s-weak} by $\varsigma \partial_t \teta_s$, the time
derivative of \eqref{eqIa} by $\varsigma \uu_t$, and the time
derivative of $\eqref{eqIIa}$ by $\varsigma \chi_t$. However, to
keep calculations simpler we have postponed this procedure to the
rigorous proof of Proposition~\ref{prop:2.1} in Sec.~\ref{rigoroso}.
\end{remark}
\subsection{Proof of Theorem~\ref{th:2.1}}
\label{ss:3.2}
 In view of estimates~\eqref{aprio-long3} and
 \eqref{est:1},  for all $\delta >0$ the trajectory
$\{ (\teta(t), \teta_s(t), \uu(t), \chi(t)), t \geq \delta\}$ of
every approximable solution $(\teta, \teta_s, \uu, \chi)$,  starting
from a quadruple  of initial data $(\teta_0,\teta_s^0,
\uu_0,\chi_0)$ as in~\eqref{cond-teta-zero}--\eqref{cond-chi-zero},
 is bounded in $V \times \hunoc \times \bfw \times H^2 (\Gamma_c)$.  Hence,
 it is relatively compact
 in $H^{1-\epsilon}(\Omega) \times H^{1-\epsilon}(\Gamma_c) \times   (H^{1-\epsilon}(\Omega))^3 \times
 H^{2-\epsilon}(\Gamma_c)$
 for all $\epsilon>0$, and with a standard argument one
 finds that the set $\omega(\teta,\teta_s, \uu,\chi)$ is a non-empty
 and compact subset of the latter product space. Thanks to~\eqref{e:imply},
$\omega(\teta,\teta_s, \uu,\chi)$ is connected  in
$H^{1-\epsilon}(\Omega) \times H^{1-\epsilon}(\Gamma_c) \times
(H^{1-\epsilon}(\Omega))^3 \times
 H^{2-\epsilon}(\Gamma_c)$
 for all $\epsilon>0$
as well, by  a well-known result in the  theory of dynamical
systems, see e.g. \cite{Haraux91}.

Now, let us fix $(\teta_\infty, \tetasinf, \uu_\infty, \chi_\infty)
\in \omega(\teta,\teta_s, \uu, \chi)$ and  an increasing sequence
$\{ t_n \} \subset (0,+\infty)$ such that $t_n \nearrow +\infty$ as
$n \up \infty$ and
 \begin{equation}
 \label{e:sense}
 \begin{aligned}
 (\teta(t_n), \teta_s (t_n),& \uu(t_n),
\chi(t_n))  \to (\teta_\infty, \tetasinf, \uu_\infty, \chi_\infty)
\\ & \text{in $ H^{1-\epsilon}(\Omega) \times H^{1-\epsilon}(\Gamma_c)
\times \left(H^{1-\epsilon}(\Omega)\right)^3 \times
H^{2-\epsilon}(\Gamma_c) $.} \end{aligned}
 \end{equation}
 Following
a well-established  procedure, for all $n \in \N$ and $T>0$ we
introduce the translated functions
\[
\begin{aligned}
\teta_n(t) := \teta(t+t_n), \quad \tetasn(t):=\teta_s (t+t_n), \quad
 \uu_n (t):= \uu(t+t_n), \quad  \chi_n(t):= \chi(t+t_n)\,,
\end{aligned}
\]
for $t \in [0,T]$. We also set for a.e. $t \in (0,T)$
\[
 \eeta_n (t):= \eeta(t+t_n),\ \
\xi_n (t):= \xi(t+t_n), \ \  \zeta_n (t):= \zeta (t+t_n), \ \
\mathbf{F}_n (t):= \mathbf{F}(t+t_n), \ \  h_n(t):= h(t+t_n).
\]
Clearly, $\eeta_n (t) \in \alpha (\uu_n (t))$ for a.e. $t \in
(0,T)$, and $\xi_n (x,t) \in \beta(\chi_n (x,t))$,
  $\zeta_n (x,t) \in \HH(\chi_n (x,t))$ for a.e. $(x,t) \in \Gamma_c \times (0,T)$.
   For later convenience, we point out that, in view
of~\eqref{cond_F-infty} and~\eqref{cond-h-infty}, for all $n \in \N$
and $T>0$
\begin{align}
& \label{e:stimafn}
 \|\mathbf{F}_n  \|_{L^\infty (0,T;W')} \leq  \|\mathbf{F} \|_{L^\infty
 (0,+\infty;W')},
 \\
  &
  \label{e:stimah_n}
  \|h_n  \|_{L^\infty (0,T;V')} \leq  \|h
  \|_{L^\infty(0,+\infty;V')}, \qquad
 \|h_n  \|_{L^1 (0,T;H)} \leq  \|h \|_{L^1(0,+\infty;H)}\,.
\end{align}
Furthermore, \eqref{e:useful-consequence} leads to
\begin{equation}
\label{e:trick} \int_0^T \| h_n \|_{V'}^2 =\int_{t_n}^{t_n+T} \| h
\|_{V'}^2 \leq \int_{t_n}^{+\infty} \| h \|_{V'}^2 \to 0 \quad
\text{as $n \to \infty$}
\end{equation}
whereby
\begin{equation}
\label{convhzero} h_n \to 0 \ \  \text{in $L^2 (0,T;V')$,} \quad h_n
\weaksto \, 0 \qquad \text{in $L^\infty (0,T;V')$,}
\end{equation}
whereas, using~\eqref{effe_infinito}, we verify that
 \begin{equation}
 \label{conve-effe-enne}
\exists\,\mathbf{F}_\infty\in \bfw'\,: \ \ \mathbf{F}_n (t) \to
\mathbf{F}_\infty \ \ \text{in $L^p (0,T;\bfw')$ for all $1\leq p
<\infty$}\,.
 \end{equation}
 Now, $(\teta_n,\tetasn, \uu_n, \chi_n, \eeta_n, \xi_n,\zeta_n)$ is
a solution to Problem~$\PP$ (with data $h_n$ and $\mathbf{F}_n$ in
place of $h$ and $\mathbf{F}$) on the interval $(0,T)$, and,
 by construction, the functions $\teta_n$, $\tetasn$, $\uu_n$ and $\chi_n$ comply
 the initial conditions
 \begin{equation}
\label{e:4.2.bis} \teta_n(0) =\teta(t_n), \ \
 \uu_n (0)= \uu(t_n) \quad \text{in $\Omega$},
\qquad \tetasn(0)= \teta_s(t_n), \ \  \chi_n(0)=\chi(t_n) \quad
\text{in $\Gamma_c$.}
\end{equation}
Exploiting~\eqref{e:sense}, we shall
 prove that the quadruple  $(\teta_\infty,\tetasinf,
\uu_\infty,\chi_\infty)$ fulfils the stationary
system~\eqref{e:stationary-system} by passing to the limit  as $n
\to \infty $ in the abovementioned initial-boundary problem.

To this aim, in view of estimates~\eqref{aprio-long}
and~\eqref{est-delta} (indeed, we may suppose without loss of
generality that, e.g.,  $t_n \geq 1$ and choose $\rho=\frac{12}{7}$
in~\eqref{est:4}), we remark that there exists a constant $K_4>0$,
independent of $n \in \N$ and $T>0$,  such that
\begin{equation}
\label{e:kappauno}
\begin{aligned}
& \| \teta_n \|_{L^\infty (0,T;V)} + \|\nabla\teta_n\|_{L^2(0,T;H)}+\| \partial_t
\ln(\teta_n)\|_{L^2 (0,T; V')}+
 \| \partial_t \teta_n \|_{L^2(0,T;L^{12/7}(\Omega))} \leq K_4\,,
 \\
& \| \tetasn \|_{L^\infty (0,T;\hunoc)} +
\|\nabla\tetasn\|_{L^2(0,T;L^2(\Gamma_c))} \leq K_4\,,
\\
&
 \| \partial_t
\ln(\tetasn)\|_{L^2 (0,T; \hunoc')}+ \| \partial_t \tetasn
\|_{L^2(0,T;L^{12/7}(\Gamma_c))} \leq K_4\,,
\\
&\|\teta_n-\tetasn\|_{L^2(0,T;L^2(\Gamma_c))}\leq K_4\,,\\
&
 \|\uu _n\|_{W^{1,\infty}(0,T;\bfw)}+ \| \partial_t \uu_n \|_{L^{2}(0,T;\bfw)}+  \|
 \widehat{\alpha}(\uu_n)\|_{L^\infty
(0,T)} \leq K_4\,,
\\
&
\begin{aligned}
 \|\chi_n   \|_{L^\infty (0,T;H^2 (\Gamma_c)) \cap W^{1,\infty}(0,T;L^2
(\Gamma_c)} &  + \| \partial_t \chi_n \|_{L^2 (0,T;H^1(\Gamma_c))}
+ \| \widehat{\beta}(\chi_n)\|_{L^\infty (0,T;L^1 (\Gamma_c))}  \\
& + \| \xi_n \|_{L^\infty (0,T;L^2 (\Gamma_c))} \leq K_4\,.
\end{aligned}
\end{aligned}
\end{equation}
Moreover, since   $0 \leq \zeta_n \leq 1$ a.e. in $\Gamma_c \times
(0,T)$ there holds
\begin{equation}
\label{e:kappadue} \| \zeta_n \|_{L^\infty (\Gamma_c \times (0,T))}
\leq 1\,.
\end{equation}
Finally, let us prove a further estimate independent of $n \in \N$
but depending on $T>0$. Exploiting~\eqref{e:kappauno} (which also
yields a bound, independent of $n\in \N$ and $T>0$,  for $\textstyle
\|\chi_n^+ \uu_n\|_{L^\infty (0,T;H^{-1/2}(\Gamma_c))}$)
and~\eqref{e:stimafn}, we argue by comparison in \eqref{eqIa}  and
conclude that
\begin{equation}
\label{e:4.2.4} \| \eeta_n \|_{L^2 (0,T;H^{-1/2}(\Gamma_c))} \leq
K_5 (T) \quad  \text{for all $n \in \N.$}
\end{equation}

Estimates~\eqref{e:kappauno}-- \eqref{e:4.2.4}, joint with standard
weak-compactness arguments, the compactness results \cite[Thm.~4,
Cor.~5]{Simon87}, and the Ascoli-Arzel\`a theorem in the framework
of the weak topologies of $\bfw$ and $H^1(\Gamma_c)$,
 yield  that there exist a subsequence  of
 $\{(\teta_n, \tetasn, \uu_n, \chi_n, \eeta_n ,\xi_n,\zeta_n) \}$
  (which we do not
relabel), and functions $(\bar{\teta}, \bar{\teta}_s,  \bar{\uu},
\bar{\chi}, \bar{\eeta}, \bar{\xi},\bar{\zeta})$ for which the
following convergences hold as $n \up \infty$
\begin{align}
\label{convtetan}
 &
\begin{aligned}
& \teta_n\weaksto\,\bar{\teta}\ \  \text{in }L^\infty(0,T;V) \cap
H^1 (0,T; L^{12/7} (\Omega)), \\ &  \teta_n \to \bar{\teta} \ \
\text{in }\mathrm{C}^0 ([0,T];H^{1-\rho}(\Omega)) \ \text{for all
}\rho \in (0,1),
\end{aligned}
 \\
 &
\label{convtetasn}
\begin{aligned}
 & \tetasn\weaksto\,\bar{\teta}_s\quad\text{in
}L^\infty(0,T;H^1 (\Gamma_c)) \cap H^1 (0,T; L^{12/7} (\Gamma_c)), \\
&  \tetasn \to \bar{\teta}_s \ \  \text{in }\mathrm{C}^0
([0,T];H^{1-\rho}(\Gamma_c)) \ \text{for all }\rho \in (0,1),
\end{aligned}
\\
& \label{convuen}
\begin{aligned}
&
\uu_n\weakto\bar{\uu}\quad\text{in }H^1(0,T;\bfw)\\
&\uu_n\to\bar{\uu}\quad\text{in }\mathrm{C}^0
([0,T];H^{1-\rho}(\Omega)^3) \quad \text{for all $\rho \in
(0,1)$\,,}
\end{aligned}
\\
& \label{convetan} \eeta_n \weakto \bar{\eeta} \quad \text{in
}L^2(0,T; (H^{-1/2}(\Gamma_c))^3)\,,
\\
& \label{convchien}
\begin{aligned} &\chi_n\weaksto\,\bar\chi\quad\text{in
}L^\infty(0,T;H^2 (\Gamma_c)) \cap H^1(0,T;\hunoc)\cap
W^{1,\infty}(0,T;L^2(\Gamma_c)),
\\
&\chi_n\rightarrow\bar{\chi}\quad\text{in }\mathrm{C}^0
([0,T];H^{2-\rho}(\Gamma_c))  \quad \text{for all $\rho \in (0,2)$},
\end{aligned}
\\
 \label{convxien}
&\xi_n\weaksto\,\bar\xi\quad\text{in }L^\infty(0,T;L^2(\Gamma_c))\,,
\\
 \label{convzetan}
&\zeta_n\weaksto\,\bar\zeta\quad\text{in
}L^\infty(0,T;L^p(\Gamma_c))\quad \text{for all $1 \leq p
<\infty$}\,.
\end{align}
Clearly, from the positivity of the sequences $\{\teta_n\}$ and
$\{\tetasn\}$ we deduce that
\[
\bar{\teta}(x,t) \geq 0 \qquad \forae\, (x,t) \in \Omega \times
(0,T), \qquad \bar{\teta}_s(x,t) \geq 0 \qquad \forae\, (x,t) \in
\Gamma_c \times (0,T)\,.
\]

 Furthermore, arguing in the same way as for~\eqref{e:trick}--\eqref{convhzero}, one sees that estimate~\eqref{est:2}
 for $\|
\partial_t \ln(\teta)\|_{L^2 (0,+\infty;V')}$  and $\| \partial_t\ln(\tetasn) \|_{L^2(0,+\infty;\hunoc')}$
 entails
 that, as $n \to \infty$,
  \begin{equation}
  \label{limln}\partial_t\ln(\teta_n) \to 0  \quad \text{in $L^2
  (0,T;V'),$} \qquad \partial_t\ln(\tetasn) \to 0  \quad \text{in $L^2
  (0,T;\hunoc')$.}
\end{equation}
Analogously, \eqref{aprio-long3} implies that
 \begin{equation}
  \label{limunt}\partial_t\uu_n \to 0  \quad \text{in $L^2
  (0,T;{\bf W}),$}
\end{equation}
and
\begin{equation}
  \label{limchint}\partial_t\chi_n \to 0  \quad \text{in $L^2
  (0,T;L^2(\Gamma_c))$}.
\end{equation}
In the same
way,  in view of \eqref{est:4} for $\teta_t$ and $\partial_t \teta_s$, we
conclude that
\begin{equation}
\label{e:precisa}
 \partial_t \teta_n \to 0 \quad \text{in $L^2 (0,T;
L^{12/7}(\Omega))$,}\qquad
\partial_t \tetasn \to 0 \quad \text{in $L^2 (0,T;
L^{12/7}(\Gamma_c))$}.
\end{equation}
Finally, from \eqref{aprio-long1} and \eqref{aprio-long2-bis}, we infer that
\begin{align}\label{convdisstetan}
& (\teta_n-\tetasn) \to 0 \quad \text{in $L^2 (0,T; L^2
(\Gamma_c))$}\,\\\no &\nabla\teta_n \to 0 \quad \text{in $L^2(0,T;
H)$},\quad\nabla\tetasn\to 0\quad\text{ in
$L^2(0,T;L^2(\Gamma_c))$}.
\end{align}
With \eqref{limunt} we get $\bar{\uu}_t =\mathbf{0}$ a.e. in
$\Omega$, so that $\bar{\uu}$ is constant in time in $\Omega$.
Hence,
\[
\bar{\uu}(x,t)=  \bar{\uu}(x,0) = \lim_{n \up \infty} \uu(x,t_n) =
\uu_\infty(x) \quad \forae\ x \in \Omega \ \ \text{for all $ t \in
[0,T]$}
\]
and analogously, from \eqref{limchint} and \eqref{e:precisa} we
deduce for all $t \in [0,T]$,
\[
\bar{\teta}(x,t)= \teta_\infty (x) \ \ \forae \  x \in \Omega\,,
\qquad \bar{\teta}_s(x,t)= \tetasinf (x)\,, \ \ \ \bar{\chi}(x,t)=
\chi_\infty (x) \ \ \forae \ x \in \Gamma_c\,.
\]
Finally, owing to \eqref{convdisstetan} we get that $\teta_\infty$
and $\tetasinf$ are constant (as
$\nabla\teta_\infty=\nabla\tetasinf={\bf 0}$) and
\begin{equation}
\label{interesting} \teta_\infty \traccia= \tetasinf \qquad \aein \
\Gamma_c\,.
\end{equation}
We also point out  that~\eqref{convuen} and~\eqref{convzetan} yield,
as $n \to \infty$,
\begin{equation}
\label{prodn} \zeta_n |\uu_n|^2 \weaksto\, \bar{\zeta} |\bar{\uu}|^2
\quad \text{in $L^\infty
  (0,T;L^{2-\rho}(\Gamma_c))$ for all $0<\rho<2$.}
\end{equation}
 Further, combining~\eqref{e:precisa}
with~\eqref{e:kappauno} and~\eqref{hyp-lambda}, we easily conclude
that $\partial_t \lambda(\chi_n)= \lambda'(\chi_n)
\partial_t \chi_n \to 0 $ in $L^2 (0,T; L^{2-\rho}(\Gamma_c))$ for
all $\rho \in (0,2)$. Arguing in the same way as  in the proof
of~\cite[Prop.~4.4]{bbr3}, we
use~\eqref{convtetan}--\eqref{convchien}
and~\eqref{limln}--\eqref{prodn} to pass to the limit in
Problem~$\PP$ and, taking into
account~\eqref{convhzero}--\eqref{conve-effe-enne}, we  conclude
that, almost everywhere in $(0,T)$, the functions $(\teta_\infty,
\tetasinf,  \uu_\infty, \chi_\infty, \bar{\eeta},
\bar{\xi},\bar{\zeta})$ fulfil \eqref{cuore} and
 \begin{align}
 &
 \label{equ-inft}
a(\uu_\infty,\vv)+ \int_{\Omega} \teta_\infty \dive (\vv) +
\int_{\Gamma_c}(\chi_\infty^+\uu_\infty +\bar\eeta )  \cdot{\bf
v}=\mathbf{F}_\infty
 \ \  \forall \vv\in \bfw\,,
\\
\label{eq-chi-inft}
 &-\Delta\chi_\infty+ \bar\xi
+\sig'(\chi_\infty)=-\lambda'(\chi_\infty)\tetasinf-\frac 1 2
\bar{\zeta}\vert\uu_\infty\vert^2 \quad\hbox{a.e. in } \Gamma_c\,,
\end{align}
joint with the no-flux boundary conditions~\eqref{bordo1}. Combining
convergences~\eqref{convchien}, \eqref{convxien},
and~\eqref{convzetan} with the strong-weak closedness of the graphs
$\beta$ and $\HH$, we find that
\begin{equation}
\label{inclu-xi} \bar{\xi} \in \beta(\chi_\infty), \quad \bar{\zeta}
\in \HH(\chi_\infty)
 \qquad
\aein\,\Gamma_c \times (0,T)\,,
\end{equation}
 while, in view of the maximal monotonicity
of the operator (induced by) $\alpha$ on $L^2
(0,T;(H^{-1/2}(\Gamma_c))^3)$ and of~\cite[Lemma~1.3, p.~42]{barbu},
we deduce that
\begin{equation}
\label{inclu-eta}
 \bar{\eeta} \in \alpha(\uu_\infty) \quad
\text{in $(H^{-1/2}(\Gamma_c))^3$} \ \ \aein\, (0,T)
\end{equation}
from the inequality
$$
\limsup_{n \up \infty}\int_0^t \int_{\Gamma_c} \eeta_n \cdot \uu_n
\leq \int_0^t \int_{\Gamma_c} \bar\eeta \cdot \uu_\infty \quad
\text{for all }t \in (0,T)\,,
$$
 which can be verified in the same way as in  the proof
of~\cite[Prop.~4.4]{bbr3}.

\begin{remark}
\upshape Notice that
$(\teta_\infty,\tetasinf,\uu_\infty,\chi_\infty)$ solve the
stationary problem \eqref{cuore}-\eqref{e:stat2} associated with the
evolution system \eqref{teta-weak}-\eqref{bordo1}. Moreover, we
observe that \eqref{limunt}--\eqref{convdisstetan} entail that
dissipation vanishes in the limit  as $t \to \infty$ (see
\eqref{pseudoO} and \eqref{pseudoG}).
\end{remark}
 \fin
\paragraph{Proof of Corollary~\ref{coro:particular-case}.}\, For
the sake of completeness, here we repeat the same argument developed
in the proof of~\cite[Prop.~2.5]{bbr2}. We test the first
of~\eqref{e:stat2} by $(\chi_\infty - m_*)$ and integrate on
$\Gamma_c$.  We   obtain
\[
\begin{aligned}
\int_{\Gamma_c} |\nabla(\chi_\infty - m_*)|^2 &  +
\int_{\Gamma_c}\xi_\infty (\chi_\infty - m_*) \\ &  =
-\int_{\Gamma_c}\left( \sigma'(\chi_\infty) +
\lambda'(\chi_\infty)\tetasinf + \frac12 \zeta_\infty
|\uu_\infty|^2\right)(\chi_\infty - m_*) \leq 0\,,
\end{aligned}
\]
the latter inequality due to~\eqref{e:spec2}--\eqref{e:spec3} and
the fact that  $m_* \leq \chi_\infty \leq m^*$ a.e. in $\Gamma_c$.
On the other hand, the second term on the left-hand side of the
above inequality is non-negative by monotonicity, so that we deduce
that $\nabla(\chi_\infty - m_*) \equiv 0$ a.e. in $\Gamma_c$. Thus,
there exists some constant $\varrho \geq m_*$ such that $\chi_\infty
\equiv \varrho$  \, a.e. in $\Gamma_c$. Now,
integrating~\eqref{e:stat2} and again
recalling~\eqref{e:spec2}--\eqref{e:spec3}, we find that
\[
\int_{\Gamma_c}\xi_\infty <0\,.
\]
Hence, necessarily $\varrho=m_*$, and~\eqref{chii-infi-det} ensues.
 \fin

\section{Rigorous estimates}\label{rigoroso}

This section is devoted to the rigorous proof of
Theorem~\ref{th:1.1} and of Proposition~\ref{prop:2.1}. Thus, in
Sec.~\ref{ss:a.1} we shall specify the variational problem
(depending on two  parameters $\eps>0$ and $\mu>0$) approximating
Problem~$\PP$ and outline the main steps of the proof of its global
well-posedness. Hence, in Sec.~\ref{ss:a.2} we shall prove some
estimates on the approximate solutions, which are independent of the
parameters $\eps$ and $\mu$ and in fact hold on $(0,+\infty)$. This
shall enable us to pass to the limit in the approximate problem
first as $\eps \searrow 0$, and secondly as $\mu \searrow 0$, and to
conclude the proof of Theorem~\ref{th:1.1} in Sec.~\ref{ss:a.3}, and
of Proposition~\ref{prop:2.1} in Sec.~\ref{ss:a.4}.

  Since the (double) approximation procedure for Problem~$\PP$ strongly
relies on the usage of Yosida regularizations of the nonlinear
operators $\ln$, $\alpha$~,  $\beta$, and of the Heaviside operator
$\HH$,  in the following section we recapitulate some related
preparatory results.
\subsection{Recaps on Yosida regularizations}
\label{ss:a0}
\begin{description}
\item[\textbf{Regularization of $\ln$.}]
For fixed $\mu>0$, we denote by
\begin{equation}
\label{e:max-mon-1} \mathsf{r}_{\mu}:= (\Id + \mu \ln)^ {-1}: \R \to
\R
\end{equation}
 the resolvent operator
associated with the logarithm $\ln $ (where $\Id: \R \to \R$ is the
identity function), and  recall that $ \mathsf{r}_{\mu}: \R \to
(0,+\infty)$ is a contraction. With easy calculations, one sees that
 $\mathsf{r}_{\mu}(1) =1$, so that, by contractivity, there holds
 \begin{equation}
 \label{e:interesting}
\mathsf{r}_{\mu}(x) \leq |x| + 2 \qquad \text{for all $x \in \R$.}
 \end{equation}
 The Yosida regularization of $\ln$ is then  defined
by
\begin{equation}
\label{e:max-mon-2} \calelleveps:=\frac{\Id -\mathsf{r}_{\mu}}{\mu}:
\R \to \R\,.
\end{equation}
It follows from~\cite[Prop.~2.6]{brezis73} that for all $\mu>0$ the
function  $\calelleveps: \R \to \R$ is non-decreasing and Lipschitz
continuous, with Lipschitz constant $1/\mu$.

For later convenience, as in~\cite{bcfg1} we also introduce the
following function
\begin{equation}
\label{def-imu} \Imu (x) := \int_0^x s \calelleveps'(s) \, \dd s\,.
\end{equation}
We point out that, since $\Imu$ is decreasing on $(-\infty,0)$ and
increasing on $(0,+\infty)$, there holds for all $\mu>0$
\begin{equation}
\label{e:immediate-to-check} \Imu (x) \geq \Imu(0)=0 \qquad \text{
for all $ x \in \R.$}
\end{equation}
 The following result collects some
properties of $\calelleveps$ and $\Imu$ which shall play a crucial
role in the proof of Proposition~\ref{prop:large-time}.
\begin{lemma}
\label{l:new-lemma} The following inequalities hold:
\begin{subequations}
\label{ine-lprimo}
\begin{align}
& \label{ine-lprimo-1} \exists\, \mu_*>0\,:  \ \ \forall\, \mu \in
(0,\mu_*) \ \ \forall\, x>0 \quad \calelleveps'(x) \leq
\frac{2}{x}\,,
\\
& \label{ine-lprimo-2}\forall\, \mu >0 \ \ \forall\, x \in \R \quad
\calelleveps'(x) \geq \frac{1}{|x| + 2 + \mu}\,.
\end{align}
\end{subequations}
As a consequence, $\Imu$ satisfies
\begin{subequations}
\label{ine-imu}
\begin{align}
& \label{ine-imu-1} \exists\, \mu_*>0\,:  \ \ \forall\, \mu \in
(0,\mu_*) \ \ \forall\, x \geq 0  \quad \Imu(x) \leq  2x\,,
\\
& \label{ine-imu-2} \exists\, C_1, \, C_2 >0 \, : \ \  \forall\, \mu
>0 \ \ \forall\, x \in \R \quad \Imu(x) \geq C_1 |x|
 - C_2\,.\end{align}
\end{subequations}
\end{lemma}
%
\noindent {\em Proof.} \,
 \textbf{Ad \eqref{ine-lprimo}.}\, Using the
 definitions~\eqref{e:max-mon-1} and \eqref{e:max-mon-2} of $\mathsf{r}_{\mu}
 $ and $\calelleveps$ and repeating the calculations
  in the proof of~\cite[Lemma~4.2]{bcfg1}, it is possible
 to show that
\[
\calelleveps' (x)  = \frac{1}{\mathsf{r}_{\mu}(x) + \mu} \qquad
\text{for all $ x \in \R $,}
 \]
 which, combined with~\eqref{e:interesting},
 yields~\eqref{ine-lprimo-2}. For the proof of~\eqref{ine-lprimo-1},
 which
 follows the very same lines, we directly refer
 to~\cite[Lemma~4.2]{bcfg1}.
\\
\textbf{Ad \eqref{ine-imu}.}\, Estimate~\eqref{ine-imu-1}  is an
immediate consequence of~\eqref{ine-lprimo-1} and of the definition
of $\Imu$. We shall now prove~\eqref{ine-imu-2} for $x \geq 0$ (the
inequality in the case $x < 0$ being  completely analogous). Indeed,
from the inequality
\[
s \calelleveps'(s)  \geq 1 - \frac{\mu + 2}{ s+ \mu +2} \qquad
\text{for all $s \geq 0$}
\]
(which is an immediate consequence of~\eqref{ine-lprimo-2}), we
deduce that
\[
\begin{aligned}
\Imu(x) & \geq x -\int_0^x \frac{\mu + 2}{s+ \mu + 2} \, \dd s
\\
 & = x - (\mu +2) \ln(x + \mu + 2) + (\mu +2) \ln( \mu + 2) \geq C_1 x - C_2\,,
\end{aligned}
\]
for some suitable positive constants $C_1$ and $C_2$.
 \fin
\item[\textbf{Regularization of $\alpha$.}]
For fixed $\mu>0$, we shall denote by
\[
\begin{aligned}
 \alpha_\mu :(H^{1/2}(\Gamma_c))^3 \to
(H^{-1/2}(\Gamma_c))^3 \  \ \text{the $\mu$-Yosida regularization of
$ \alpha$,}
\end{aligned}
\]
and recall that, by~\cite[Prop.~II.1.1]{barbu}, the operator
$\alpha_\mu: (H^{1/2}(\Gamma_c))^3 \to (H^{-1/2}(\Gamma_c))^3$ is
single-valued, monotone, bounded and demi-continuous. Being
$\widehat{\alpha}$ a non-negative functional by~\eqref{hyp:alpha},
 its primitive \[
 \widehat{\alpha}_\mu(\uu):=
\min_{\vv \in \overline{D(\alpha)}}
\left(\frac{\|\uu-\vv\|_{H^{1/2}(\Gamma_c)}^2}{2\mu} +
\widehat{\alpha}(\vv)\right) \qquad \text{for all $\uu \in
(H^{1/2}(\Gamma_c))^3$}
\]
(which is Fr\'echet differentiable on $(H^{1/2}(\Gamma_c))^3$),
fulfils
\begin{equation}
\label{a1.2} 0 \leq  \widehat{\alpha}_\mu(\uu) \leq
\widehat{\alpha}(\uu) \qquad \text{for all $\uu \in
\overline{D(\alpha)}$.}
\end{equation}
\item[\textbf{Regularization of $\beta$.}]
We shall also use
\[
 \beta_\mu :\R  \to \R \  \ \text{the $\mu$-Yosida
regularization of $ \beta$.}
\]
With straightforward calculations one verifies that, thanks
to~\eqref{e:consequences},
 the Yosida approximation
\[
\widehat{\beta}_\mu(x):= \min_{y \in \overline{D(\beta)}}
\left(\frac{|y-x|^2}{2\mu} + \widehat{\beta}(y)\right) \qquad
\text{for all $ x \in \R$}\,
\]
of $\widehat{\beta} $ verifies
\begin{equation}
\label{a1.1}
\begin{cases}
 \forall\, R>0 \ \ &\exists\, C_{R}>0 \  \  \forall\, \mu>0 \
\forall\, x \geq 0\, : \ \ \widehat{\beta}_\mu(x) \geq \frac{R}{2\mu
R + 1 } x^2 - C_R\,,
\\
& \exists\,\overline{C}>0 \ \ \forall\, \mu>0 \ \forall\, x <0\, : \
\ \widehat{\beta}_\mu(x) \geq  \frac{x^2}{2\mu} - \overline{C}\,.
\end{cases}
\end{equation}
\item[\textbf{Regularization of $\HH$.}]
For fixed $\mu>0$, we shall denote by
\[
\begin{aligned}
& \HHmu = \ppmu'  :\R  \to \R \  \ \text{the $\mu$-Yosida
regularization of $ \HH$,}
\\
&  \ppmu :\R  \to \R \  \ \text{the $\mu$-Yosida approximation of $
(\cdot)^+$.}
\end{aligned}
\]
In particular, it can be checked that
\begin{equation}
\label{formpmu} \ppmu(x) = \begin{cases} 0 \qquad & \text{if }x \leq
0,
\\
\displaystyle \frac{x^2}{2\mu} \qquad & \text{if }0 <x <\mu,
\\
\displaystyle x-\frac\mu2 \qquad & \text{if }x \geq \mu.
\end{cases}
\end{equation}
Using the definition of $\HHmu$ and $\ppmu$, it is straightforward
to verify that
\begin{align}
& \label{def-hmu} 0 \leq \HHmu(x) \leq 1 \qquad \text{for all $x \in
\R$,}
\\
& \label{defpmu} 0 \leq \ppmu(x) \leq (x)^+  \qquad \text{for all $x
\in \R$.}
\end{align}
\end{description}

\subsection{Approximation of Problem~$\PP$}
\label{ss:a.1}
\paragraph{A  double approximation procedure.} Let $\eps, \, \mu >0$
be two strictly positive parameters. We consider the approximation
of Problem~$\PP$ obtained in the following way:
\begin{enumerate}
\item
 we add to~\eqref{teta-weak} the  regularizing viscosity
term $\eps \mathcal{R} (\teta_t)$  and to~\eqref{teta-s-weak} the
viscosity term $\eps \mathcal{R}_{\Gamma_c} (\partial_t \teta_s)$
($\mathcal{R}$ and $ \mathcal{R}_{\Gamma_c}$ being the Riesz
operators introduced in Notation~\ref{not:2.1});
\item we replace the operators $\alpha$ in~\eqref{eqIa}, $\beta$ and $\HH$ in~\eqref{eqIIa} with
their Yosida regularization $\alpha_\mu
:\left(H^{1/2}(\Gamma_c)\right)^3 \to
\left(H^{-1/2}(\Gamma_c)\right)^3$, $\beta_\mu : \R \to \R$, and
$\HHmu: \R \to \R$; accordingly, we replace the term $\chi^+ \uu$ in
equation~\eqref{eqIa} by $\ppmu (\chi) \uu$;
\item
  both in~\eqref{teta-weak} and in~\eqref{teta-s-weak}  we replace the
  logarithm
$\ln$  with its Yosida regularization $\calelleveps$.
\end{enumerate}
\noindent
\paragraph{Approximate initial data.}
In order to properly state our approximate problem, depending on the
parameters $\eps>0$ and $\mu>0$, we shall need some enhanced
regularity on the initial data for  $\teta$ and $\teta_s$. The
following result concerns the construction of sequences of suitable
approximate initial data $\{\tetazeroemu\}$ and
$\{\tetaessezeroemu\}$, which in  fact depend on the parameter
$\eps>0$ only.
\begin{lemma} \label{dati-iniziali-lemma}
 Assume that the initial
data $\teta_0$ and $\teta_s^0$ respectively comply
with~\eqref{cond-teta-zero} and~\eqref{cond-teta-esse-zero}.
 Then,
 \begin{enumerate}
\item there exists a sequence $\{\tetazeroemu\} \subset V$ such that
for all $\eps>0$ and $\mu \in (0,\mu_*)$ ($\mu_*$ being as in
Lemma~\ref{l:new-lemma}) there hold
\begin{subequations}
\label{tetazeroapp}
\begin{align}
& \label{tetazeroapp-1} \eps^{1/2} \| \tetazeroemu \|_{V} \leq
\|\teta_0\|_{V'}\,,
\\
& \label{tetazeroapp-2} \int_{\Omega}\Imu( \tetazeroemu )  \leq 2
\|\teta_0\|_{L^1 (\Omega)}\,,
\\
& \label{tetazeroapp-3}  \tetazeroemu  \to \teta_0 \quad \text{in
$L^{\bar{p}} (\Omega)$ as $\eps \searrow 0$}\,,
\\
& \label{tetazeroapp-4}  \calelleveps(\tetazeroemu)  \to
\calelleveps(\teta_0) \quad \text{in $L^{\bar{p}} (\Omega)$ as $\eps
\searrow 0$}\,;
\end{align}
\end{subequations}
\item there exists a sequence $\{\tetaessezeroemu\} \subset \hunoc$ such that
for all $\eps>0$ and $\mu \in (0,\mu_*)$
\begin{subequations}
\label{tetazeroapps}
\begin{align}
& \label{tetazeroapps-1} \eps^{1/2} \| \tetaessezeroemu \|_{\hunoc}
\leq \|\tetaessezero\|_{\hunoc'}\,,
\\
& \label{tetazeroapps-2} \int_{\Gamma_c}\Imu( \tetaessezeroemu )
\leq 2 \|\tetaessezero\|_{L^1 (\Gamma_c)}\,,
\\
& \label{tetazeroapps-3}  \tetaessezeroemu \to \tetaessezero \quad
\text{in $L^{\bar{q}} (\Gamma_c)$ as $\eps \searrow 0$}\,,
\\
& \label{tetazeroapps-4}  \calelleveps(\tetaessezeroemu)  \to
\calelleveps(\teta_s^0) \ \text{in $L^{\bar{q}}(\Gamma_c)$ as $\eps
\searrow 0$}\,.
\end{align}
\end{subequations}
 \end{enumerate}
\end{lemma}
\noindent \emph{Proof.} We shall carry out the construction of the
sequence $\{\tetazeroemu \}$, the passages for $\{ \tetaessezeroemu
\}$ being completely analogous. For all $\eps,\, \mu>0$, we let
$\tetazeroemu  \in V$ be the solution of the variational equation
\begin{equation}
\label{e:variational-eps} \int_{\Omega} \tetazeroemu v + \eps^{1/2}
\int_{\Omega} \nabla\tetazeroemu \nabla v= \int_{\Omega} \teta_0 v
\qquad \text{for all $v \in V$}
\end{equation}
(notice that $\tetazeroemu \in V$ is well-defined, thanks
to~\eqref{sto-anche-li}).
 Then, \eqref{tetazeroapp-1} can be straightforwardly proved,
 together with
 \begin{equation}
\label{oh-la-la}\tetazeroemu \to \teta_0 \qquad \text{in $V'$ as
$\eps \searrow 0$.}
 \end{equation}
Furthermore, the maximum principle shows that, being $\teta_0 >0$
a.e. in $\Omega$ thanks to the second of~\eqref{cond-teta-zero},
\[
\tetazeroemu \geq 0 \qquad \text{a.e. in $\Omega$.}
\]
Then, testing~\eqref{e:variational-eps} by $1$, one  obtains
\begin{equation}
\label{e:proof1} \|\tetazeroemu \|_{L^1 (\Omega)} = \|\teta_0
\|_{L^1 (\Omega)} \quad \text{for all $\eps>0$.}
\end{equation}
Therefore, \eqref{tetazeroapp-2} follows from
combining~\eqref{e:proof1} with~\eqref{ine-imu-1}.
 Finally, to check~\eqref{tetazeroapp-3} we test~\eqref{e:variational-eps}
by $(\tetazeroemu)^{{\bar{p}}-1} $ and with easy calculations deduce
that $\| \tetazeroemu \|_{L^{\bar{p}} (\Omega)} \leq \| \teta_0
\|_{L^{\bar{p}} (\Omega)}$. Combining this estimate
with~\eqref{oh-la-la} using the uniform convexity of $L^{\bar{p}}
(\Omega)$, we conclude the strong convergence of $\{  \tetazeroemu
\}$ in $L^{\bar{p}} (\Omega)$. Hence, \eqref{tetazeroapp-4} is a
direct consequence of~\eqref{tetazeroapp-4} and the Lipschitz
continuity of $\calelleveps$.
 \fin
\paragraph{Variational formulation of the approximate problem.}
 We thus obtain the following boundary value
problem, which we directly  state on the half-line $(0,+\infty)$ in
view of the long-time a priori estimates of
Proposition~\ref{prop:large-time}.
\medskip\par\noindent
{\bf Problem $\pepsmu$}.  Given a quadruple of  initial data
$(\tetazeroemu,\tetaessezeroemu, \uu_0, \chi_0)$, $\uu_0$ and
$\chi_0$ being as in \eqref{cond-uu-zero}--\eqref{cond-chi-zero},
and $\tetazeroemu$ and $\tetaessezeroemu$
fulfilling~\eqref{tetazeroapp} and~\eqref{tetazeroapps}
respectively,
 find  functions  $( \teta, \teta_s,
\uu,\chi)$, with the regularity
\[
\begin{aligned}
& \teta \in H^1 (0,T; V), \ \ \teta_s \in  H^1 (0,T; \hunoc), \ \
 \uu \in H^1 (0,T;\bfw)\,,
\\
 &
 \chi   \in L^{2}(0,T;H^2 (\Gamma_c))  \cap L^{\infty}(0,T;H^1
(\Gamma_c)) \cap H^1 (0,T;L^2 (\Gamma_c))
 \end{aligned}
\]
for all $T>0$,  fulfilling the initial conditions
\begin{equation}
\label{approx-ini-teta} \teta(0)= \tetazeroemu \quad \aein\, \Omega,
\qquad \qquad \teta_s (0)= \tetaessezeroemu \quad \aein\, \Gamma_c
\end{equation}
and~\eqref{iniu}--\eqref{inichi}, and the equations
\begin{align}
& \label{teta-weak-app}
\begin{aligned}
\eps \int_{\Omega} \teta_t \, v &+ \int_{\Omega}
\partial_t\calelleveps(\teta) v -\int_{\Omega} \dive(\uu_t) \, v
+\eps\int_{\Omega} \nabla \teta_t \, \nabla v +\int_{\Omega} \nabla
\teta \, \nabla v
\\ &+ \int_{\Gamma_c} k(\chi)
(\teta-\teta_s) v = \pairing{}{}{h}{v} \quad \forall\, v \in V \
\hbox{ a.e. in }\, (0,+\infty)\,,
\end{aligned}
\\ & \label{teta-s-weak-app}
\begin{aligned}
&\eps \int_{\Gamma_c} \partial_t \teta_s \, v + \int_{\Gamma_c}
\partial_t\calelleveps(\teta_s) v
+\eps\int_{\Gamma_c}\nabla\partial_t{\teta_s}\nabla v
\\ &-\int_{\Gamma_c}
\partial_t \lambda(\chi) \, v    +\int_{\Gamma_c} \nabla \teta_s  \, \nabla
v = \int_{\Gamma_c} k(\chi) (\teta-\teta_s) v  \quad \forall\, v \in
\hunoc  \ \hbox{ a.e. in }\, (0,+\infty)\,,
\end{aligned}
\\
& \label{eqIappr}
\begin{aligned}
 b({\bf {u}}_t,\vv)+a(\uu,\vv)+
\int_{\Omega} \teta \dive (\vv)& + \int_{\Gamma_c}(\ppmu(\chi){\bf
u} +\alpha_\mu (\uu) ) \cdot{\bf v} \\ &=
\pairing{\bfw'}{\bfw}{\mathbf{F}}{\vv}
 \quad \forall \vv\in \bfw \ \hbox{ a.e. in }\, (0,+\infty)\,,
 \end{aligned}
\\
 &
\label{eqIIappr}
\begin{aligned}
 \chi_t -\Delta\chi+ \beta_\mu (\chi)
+\sig'(\chi) & =-\lambda'(\chi) \teta_s-\frac 1 2\vert{\bf
u}\vert^2\HHmu(\chi) \quad\hbox{a.e. in } \Gamma_c\times
(0,+\infty)\,,
\\
 & \partial_{\nn_s}
  \chi=0\text{ a.e. in } \partial\Gamma_c\times (0,+\infty)\,.
\end{aligned}
\medskip
\end{align}
\begin{remark}\label{piufacile}
\upshape Note that the approximate system
\eqref{teta-weak-app}-\eqref{eqIIappr} presents fewer technical
difficulties than the analogous approximate version introduced in
\cite{bbr3}. This is mainly due to the fact that, as we deal with
the specific choice of $\ln(\teta)$ and $\ln(\teta_s)$ in the
entropy  (cf. Remark~\ref{scelgolog}),  we are allowed to directly
introduce the Yosida regularization of the logarithm, instead of the
more intricate approximating procedure exploited in \cite{bbr3}.
\end{remark}
\subsubsection{Outline of the proof of global well-posedness for
Problem~$\pepsmu$.}\, Since Problem~$\pepsmu$ only slightly differs
from the approximate problem considered in~\cite{bbr3},  the global
well-posedness  for Problem~$\pepsmu$, on any interval $(0,T)$, can
be obtained arguing in the very same way as in~\cite[Sec.~3]{bbr3},
to which we refer the reader for all details. Here, we shall just
sketch the main steps of the proof.
\begin{description}
\item[\textbf{Step $1.$}]
 First of all, one proves the existence of a local (in time)
solution to (the Cauchy problem for) Problem~$\pepsmu$ with the use
of a Schauder fixed point argument. This involves establishing
intermediate well-posedness results for each of the approximate
equations, for which we refer to the calculations developed
in~\cite[Sec.~3.2]{bbr3}, and defining a solution operator which
complies with the conditions of the Schauder fixed point theorem
(cf.~\cite[Sec.~3.3]{bbr3}).
\item[\textbf{Step $2.$}]
Next, to extend the local solution to the whole interval $(0,T)$,
one needs global (in time) a priori estimates. The latter
substantially coincide with the ones formally performed in the proof
of Proposition~\ref{prop:2.1} and shall be repeated on the
approximate system~\eqref{teta-weak-app}--\eqref{eqIIappr} within
the proof of Proposition~\ref{prop:large-time}.
\item[\textbf{Step $3.$}]
Finally, uniqueness of solutions to Problem~$\pepsmu$ follows from
the very same contraction estimates performed
in~\cite[Sec.~3.5]{bbr3}.
\end{description}
\begin{remark}
\upshape \label{small-explanation} We  briefly justify our
construction of the approximate Problem~$\pepsmu$,
  referring
to~\cite[Rems.~3.2, 4.2]{bbr3}
 for  more comments.

The viscosity terms $\eps \mathcal{R}(\teta_t)$ and  $\eps
\mathcal{R}\traccia (\partial_t \teta_s)$ have been inserted
in~\eqref{teta-weak-app} and~\eqref{teta-s-weak-app}, respectively,
for technical reasons,   related to the fixed point construction of
a local solution for Problem~$\pepsmu$.  In particular, the
contributions $\eps \mathcal{R}(\teta_t)$ and  $\eps
\mathcal{R}\traccia (\partial_t \teta_s)$ are essential to prove
uniqueness of solutions to (the Cauchy problems for) approximate
equations~\eqref{teta-weak-app} and~\eqref{teta-s-weak-app}.
Furthermore,  they also play a crucial role to make the estimate
leading to the further regularity~\eqref{est:1}--\eqref{est:4}
rigorous, see the ensuing proof of
Proposition~\ref{prop:large-time}. For the same reason, we have
replaced the maximal monotone operators $\ln$, $\alpha$ and $\beta$
by their Yosida regularizations, and correspondingly substituted the
coupling terms $\chi^+ \uu$ and $-1/2 \zeta |\uu|^2 $ in
equations~\eqref{eqIa} and~\eqref{eqIIa}, with $\ppmu(\chi) \uu$ and
$-1/2 \HHmu(\chi)|\uu|^2 $, respectively.

 Now, as in the approximation of the system considered in~\cite{bbr3},
we shall keep
  the viscosity parameter $\eps$ distinct from the Yosida parameter
 $\mu$ in both approximate equations \eqref{teta-weak-app} and
 \eqref{teta-s-weak-app}. Thus, we shall prove
the existence of solutions to Problem~$\PP$  by passing to the limit
in Problem~$\pepsmu$ first as $\eps \searrow 0$ for $\mu>0$ fixed,
and then as $\mu \searrow 0$. This procedure shall enable us
 to recover on any interval $(0,T)$ the
$L^\infty(0,T;H)$-regularity for $\ln(\teta)$ (the $L^\infty
(0,T;L^2(\Gamma_c))$-regularity for  $\ln(\teta_s)$, respectively)
by testing  the $\mu$-approximation of \eqref{teta-weak} (of
\eqref{teta-s-weak}, respectively), by the term
$\calelleveps(\teta)$ ($\calelleveps(\teta_s)$, resp.), and
 obtaining some bound  independent of the approximation
parameter $\mu$. In fact, such an estimate may be performed on
equation \eqref{teta-weak-app} (on \eqref{teta-s-weak-app}, resp.)
only when
 $\eps=0$. For, if one kept $\eps>0$, one would not obtain
estimates on $\calelleveps (\teta)$ independent of the parameters
$\eps$ and $\mu$, essentially because the term $\pairing{}{}{\eps
\mathcal{R}(\teta_t)}{\calelleveps (\teta)}$ ($\pairing{}{}{\eps
\mathcal{R}\traccia(\partial_t\teta_s)}{\calelleveps (\teta_s)}$,
resp.) cannot be dealt with by monotonicity arguments.
\end{remark}

\subsection{Estimates on solutions to Problem~$\pepsmu$}
\label{ss:a.2}
\begin{proposition}
\label{prop:large-time} Assume~\eqref{A5}--\eqref{hyp-k} and
 \eqref{e:consequences}--\eqref{cond-h-infty}.
Then,
\begin{enumerate}
\item
there
 exists a constant $K_6>0$, only depending on  the functions $\lambda$, $k$, $\sigma$,
  and on the quantity $M$~\eqref{def:M},
such that for  all $\eps,\,\mu>0$ the following estimates hold for
the family of solutions $\{(\tetame,\tetasme, \ume, \chime)
\}_{\eps,\mu}$ to Problem~$\pepsmu$:
\begin{subequations}
\label{aprio-long-epsmu}
\begin{align}
& \label{aprio-long1-epsmu} \| \nabla \tetame \|_{L^2 (0,+\infty;
H)}  + \| \nabla \tetasme \|_{L^2 (0,+\infty; L^2 (\Gamma_c))} \leq
K_6\,,
\\
 &
 \label{aprio-long2-epsmubis}
\eps^{1/2} \| \tetame\|_{L^\infty (0,+\infty; V)} + \eps^{1/2} \|
\tetasme\|_{L^\infty (0,+\infty; \hunoc)}\leq K_6\,,
\\
 &
 \label{aprio-long2-epsmu}
\| \tetame \|_{L^\infty (0,+\infty; L^1 (\Omega))} + \| \tetasme
\|_{L^\infty (0,+\infty; L^1 (\Gamma_c))}\leq K_6\,,
\\
& \label{aprio-long2-bis-epsmu}
 \| \tetame -\tetasme \|_{L^2 (0,+\infty; L^2 (\Gamma_c))} \leq K_6\,,
\\
& \label{aprio-long3-epsmu}
  \| \partial_t\ume\|_{L^2 (0,+\infty;\bfw)} + \| \ume \|_{L^\infty
(0,+\infty;\bfw)} + \|  \widehat{\alpha}_\mu(\ume)\|_{L^\infty
(0,+\infty)} \leq K_6\,,
\\
& \label{aprio-long4-epsmu}
\begin{aligned}
 & \|
\chime \|_{L^\infty (0,+\infty;H^1 (\Gamma_c))}  +
\|\partial_t\chime  \|_{L^2 (0,+\infty; L^2 (\Gamma_c))} \\  & + \|
\widehat{\beta}_\mu(\chime)\|_{L^\infty (0,+\infty;L^1 (\Gamma_c))}
+ \| \HHmu(\chime)\|_{L^\infty (\Gamma_c \times (0,+\infty))} \leq
K_6\,, \end{aligned}
\\
& \label{aprio-long5-epsmu} \| \partial_t \calelleveps(\tetame)
\|_{L^2 (0,+\infty; V')} + \| \partial_t  \calelleveps(\tetasme)
\|_{L^2 (0 ,+\infty; \hunoc')} \leq K_6\,.
\end{align}
\end{subequations}
\item Furthermore, for all $\delta>0$  there exist constants $K_7(\delta), \, K_{8}(\delta,\rho) >0$,
  depending on $\delta$,   on  the functions $\lambda$, $k$, $\sigma$, and on the quantity
$M$~\eqref{def:M} ($K_{8}(\delta,\rho)$ on $\rho \in (0,2)$ as
well), but independent of $\eps>0$ and $\mu>0$, such that
 the following estimates hold
 \begin{subequations}
 \label{est-delta-me}
\begin{align}
& \label{est:1-me} \begin{aligned} & \| \tetame \|_{L^\infty
(\delta,+\infty;V)}  + \| \tetasme \|_{L^\infty
(\delta,+\infty;H^1(\Gamma_c))} \\ &  + \| \chime \|_{L^\infty
(\delta,+\infty;H^2 (\Gamma_c))}   + \| \partial_t \chime\|_{L^2
(\delta,+\infty;\hunoc) \cap  L^\infty (\delta,+\infty;L^2
(\Gamma_c))} \\ & + \| \beta_\mu (\chime) \|_{L^\infty
(\delta,+\infty;L^2 (\Gamma_c))}  + \|\ume \|_{W^{1,\infty}
(\delta,+\infty;\bfw)} \leq K_7(\delta)\,,
\end{aligned}
\\
 & \label{est:4-me}
 \|
\partial_t \tetame \|_{L^2 (\delta,+\infty;L^{12/7}(\Omega))} + \|
\partial_t \tetasme\|_{L^2 (\delta,+\infty;L^{2-\rho}(\Gamma_c))}  \leq
K_{8}(\delta,\rho) \,,
\end{align}
\end{subequations}
\end{enumerate}
\end{proposition}
\noindent \emph{Proof.}
 We shall prove~\eqref{aprio-long-epsmu}
and~\eqref{est-delta-me} by performing on Problem~$\pepsmu$ the same
a priori estimates as in the proof of Proposition~\ref{prop:2.1} and
obtaining bounds independent  of $\eps,\, \mu>0$ and of $t \in
(0,+\infty)$.
\paragraph{First  estimate.}\,
 We test~\eqref{teta-weak-app} by $ \tetame$,
\eqref{teta-s-weak-app} by $ \tetasme$, \eqref{eqIappr} by $
\partial_t\ume$ and~\eqref{eqIIappr} by $
\partial_t\chime$, add the resulting relations and integrate them on
the interval $(0,t)$, with $t \in (0,+\infty)$. Basically, we put
forth the same calculations as
throughout~\eqref{e:formal1}--\eqref{e2.1.2}, up to the following
changes. Instead of the formal identities~\eqref{e:formal1}, using
the definition~\eqref{def-imu} of $\Imu$ we find
\[
\begin{aligned}
 \int_0^t\int_{\Omega} \partial_t \calelleveps(\tetame) \tetame & =
\int_0^t\int_{\Omega} \Imu'(\tetame)  \partial_t \tetame \\ & =
 \int_{\Omega} \Imu
\left( \tetame (t) \right) - \int_{\Omega} \Imu
\left(\tetazeroemu\right) \\ & \geq C_1 \|\tetame (t)\|_{L^1
(\Omega)} - C_2 - 2 \| \teta_0 \|_{L^1 (\Omega)}\,,
\end{aligned}
\]
the last inequality ensuing from~\eqref{ine-imu-2}
and~\eqref{tetazeroapp-2}. In the same way, we have
\[\int_0^t\int_{\Gamma_c} \partial_t \calelleveps(\tetasme) \tetasme
\geq  C_1 \|\tetasme (t)\|_{L^1 (\Gamma_c)} - C_2 - 2 \|
\tetaessezero \|_{L^1 (\Gamma_c)}\,,
\]
whereas we estimate
\[
\int_0^t \langle \eps \mathcal{R}(\partial_t\tetame), \tetame
\rangle = \frac\eps2\| \tetame(t) \|_{V}^2 - \frac\eps2\|
\tetazeroemu \|_{V}^2 \geq \frac\eps2\| \tetame(t) \|_{V}^2 -\frac12
\| \teta_0\|_{V'}^2
\]
in view of~\eqref{tetazeroapp-1}, and analogously for $\tetasme$.
 Then, we perform the same
 passages as in~\eqref{e:rig1}--\eqref{e:to-be-cited}, up to replacing $\widehat{\alpha}$,
$(\cdot)^+$, and
  $\widehat{\beta}$ with their approximations
  $\widehat{\alpha}_\mu$, $\ppmu$, and
  $\widehat{\beta}_\mu$.
Repeating the very same calculations as
in~\eqref{e2.1.1}--\eqref{e2.1.2} and taking into account some
cancellations as well as~\eqref{tetazeroapp-1}
and~\eqref{tetazeroapps-1}, we arrive at
\[
\begin{aligned}
\frac\eps2 \| \tetame(t)\|_{V}^2  & +  \| \tetame(t) \|_{L^1
(\Omega)} + \frac12\int_0^t \| \nabla \tetame \|_H^2 + \int_0^t
\int_{\Gamma_c} k(\chime) (\tetame-\tetasme)^2\\ &  \frac\eps2 \|
\tetasme(t) \|_{\hunoc}^2 + \| \teta_s(t) \|_{L^1 (\Gamma_c)} +
\int_0^t \| \nabla \teta_s \|_{L^2 (\Gamma_c)}^2 \\ & + C_b \int_0^t
\|
\partial_t\ume\|_{\bfw}^2  + \frac{C_a}4 \| \ume(t) \|_{\bfw}^2 +
\widehat{\alpha}_\mu (\ume(t)) + \frac12\int_{\Gamma_c}
\ppmu(\chime(t)) |\ume(t)|^2 \\ & + \int_0^t \| \partial_t\chime
\|_{L^2 (\Gamma_c)}^2 + \frac12 \| \nabla \chime(t) \|_{L^2
(\Gamma_c)}^2 + \int_{\Gamma_c}\left( \widehat{\beta}_\mu(\chime(t))
+ \sigma(\chime(t)) \right) \\ &
\begin{aligned}
 \leq  C\Big( & 1  +  \| \teta_0\|_{L^1 (\Omega)} + \| \teta_0 \|_{V'} + \|
\teta_s^0 \|_{L^1 (\Gamma_c)} +\| \teta_s^0 \|_{\hunoc'} \\ &    +\|
\uu_0 \|_{\bfw}^2 + \widehat{\alpha} (\uu_0)  + \| \nabla \chi_0
\|_{L^2 (\Gamma_c)}^2
\\ &  +
\| \chi_0\|_{L^2 (\Gamma_c)} \| \uu_0 \|_{\bfw}^2  + \|
\widehat{\beta}(\chi_0) \|_{L^1 (\Gamma_c)}+ \| \sigma(\chi_0)
\|_{L^1 (\Gamma_c)} +   \int_0^t  \| h \|_{V'}^2 \\ & + \|
\mathbf{F}\|_{L^\infty (0,+\infty;\bfw')}^2\Big)
  + \frac1{|\Omega|^{1/2}} \int_0^t \| h \|_{V'} \|
\teta\|_{L^1 ((\Omega)}  + \int_0^t \|\mathbf{F}_t\|_{\bfw'} \|\ume
\|_{\bfw}\,.
\end{aligned}
\end{aligned}
\]
As in the proof of Proposition~\ref{prop:2.1}, from this inequality
we deduce the
bounds~\eqref{aprio-long1-epsmu}--\eqref{aprio-long2-epsmu} and
\eqref{aprio-long3-epsmu}--\eqref{aprio-long4-epsmu} (note that
estimate \eqref{beta-sigma} holds on this approximate level, too,
thanks to~\eqref{a1.1}, and the estimate for $\HHmu (\chime)$
trivially follows from~\eqref{def-hmu}). Further,
\eqref{aprio-long2-bis-epsmu} ensues from the bound
\begin{equation}
\label{e:utile-dopo} \| (k(\chime))^{1/2} (\tetame-\tetasme)\|_{L^2
(0,+\infty; L^2 (\Gamma_c))} \leq C
\end{equation}
and from~\eqref{hyp-k}. Again, the bound~\eqref{aprio-long5-epsmu}
for $\partial_t \calelleveps (\tetame)$ and $\partial_t \calelleveps
(\tetasme)$ is a consequence of  the previous estimates and a
comparison in equations~\eqref{teta-weak-app}
and~\eqref{teta-s-weak-app}.
\paragraph{Second estimate.}\,
  We test~\eqref{teta-weak-app} by $ \tanh(\cdot) \partial_t \tetame$,
  \eqref{teta-s-weak-app} by $\tanh(\cdot)\partial_t \tetasme$,
differentiate~\eqref{eqIappr} w.r.t. time and test it by
$\tanh(\cdot)\partial_t\ume$, and differentiate~\eqref{eqIIappr} and
multiply it by $\tanh(\cdot)\partial_t\chime$. Notice that such an
estimate is rigorous on this level, for $\partial_t \tetame \in V$
and $\partial_t \tetasme \in \hunoc$. Further, we sum the resulting
relations and integrate in time.
 In place of using the formal
identities~\eqref{e:formal2}, we have the following inequalities,
which are direct consequences of~\eqref{ine-lprimo-2}:
\begin{equation}
\label{pr0}
\begin{aligned} \int_0^t     & \int_{\Omega} \partial_t \calelleveps(\tetame(\cdot,r))
\tanh(r)  \partial_t \tetame(\cdot,r) \dd r \\ &  = \int_0^t
\int_{\Omega} \tanh(r) |\partial_t\tetame(\cdot,r)|^2 \calelleveps
'(\tetame(\cdot,r)) \dd r\\ &  \geq \int_0^t \int_{\Omega}\tanh(r)
\frac{ |\partial_t\tetame(\cdot,r)|^2}{|\tetame (\cdot,r)| + \mu +
2} \dd r
\\ & =4\int_0^t \int_{\Omega} \tanh(r)|\partial_t
\Thetame^{1/2} (\cdot,r)|^2 \dd r\,,
\end{aligned}
\end{equation}
and, analogously,
\begin{equation}
\label{pr1}
\begin{aligned} \int_0^t &  \int_{\Gamma_c} \partial_t \calelleveps(\tetasme(\cdot,r))
\tanh(r)  \partial_t \tetasme(\cdot,r) \dd r \geq  \int_0^t
\int_{\Gamma_c} \tanh(r)|\partial_t  \Thetasme^{1/2}(\cdot,r) |^2
\dd r\,,
\end{aligned}
\end{equation}
where we have set for all $\eps,\, \mu>0$
\begin{equation}
\label{teta-big} \begin{aligned} &  \Thetame: =|\tetame|+ \mu + 2\,,
\qquad
 \Thetasme:=|\tetasme|+ \mu + 2 \,.
\end{aligned}
\end{equation}
 Then, we observe that the remaining computations are not
affected by the factor $\tanh(\cdot)$ in a substantial way. In fact,
the same calculations as in the proof of the second (formal)
estimate in Proposition~\ref{prop:2.1} go through, up to dealing
with the integration by parts more carefully due to the presence of
$\tanh(\cdot)$. In particular, we have
\begin{equation}
\label{pr2}
\begin{aligned}
\int_0^t \int_{\Omega} & \nabla \tetame(\cdot,r) \nabla \left(
\tanh(r) \partial_t \tetame(\cdot,r) \right) \dd r \\ &  = \frac12
\int_{\Omega} \tanh(t) |\nabla \tetame(t)|^2 - \frac12 \int_0^t
\tanh'(r) \int_{\Omega} |\nabla \tetame(\cdot,r)|^2  \dd r
\\ & \geq \frac{\tanh(t)}2
\int_{\Omega} |\nabla \tetame(t)|^2 - \frac{K_6^2}2 \,,
\end{aligned}
\end{equation}
in which we have used the fact that $\tanh(0)=0$,  that $0<\tanh'(r)
\leq 1$ for all $r \in \R$, and estimate~\eqref{aprio-long1-epsmu}.
We handle the corresponding term for $\tetasme$ in the same way.
Next, we estimate
\begin{equation}
\label{pr3}
\begin{aligned}
 \int_0^t   \int_{\Gamma_c}  &   k(\chime \cdr) \tanh(r)
(\tetame \cdr - \tetasme
\cdr) \partial_t \left( \tetame \cdr - \tetasme \cdr \right) \dd r\\
\\ & = \frac{\tanh(t)}2 \int_{\Gamma_c} k(\chime (t)) \left(\tetame(t) -
\tetasme(t) \right)^2 + I_{10} + I_{11}\,,
\end{aligned}
\end{equation}
where, integrating by parts
\begin{align}
& \label{pr4}
 I_{10} = -\frac12\int_0^t \int_{\Gamma_c} \tanh'(r) k(\chime
\cdr) \left( \tetame \cdr - \tetasme \cdr \right)^2 \geq -C
\\
& \label{pr5}
\begin{aligned}
 I_{11} &  = -\frac12\int_0^t \int_{\Gamma_c} \tanh(r)
k'(\chime \cdr) \partial_t \chime \cdr \left( \tetame \cdr -
\tetasme \cdr \right)^2  \dd r\\ & \geq - \frac{L_k}{2} \int_0^t
\|\partial_t \chime  \|_{L^4 (\Gamma_c)} \|\tetame  - \tetasme
\|_{L^4 (\Gamma_c)} \|\tetame  - \tetasme \|_{L^2 (\Gamma_c)}
\\ & \geq -\nu \int_0^t \|\partial_t \chime  \|_{L^4
(\Gamma_c)}^2   - C_\nu \int_0^t \|\tetame  - \tetasme \|_{L^4
(\Gamma_c)}^2 \|\tetame  - \tetasme \|_{L^2 (\Gamma_c)}^2
\end{aligned}
\end{align}
the inequality in~\eqref{pr4} due to~\eqref{e:utile-dopo}, while the
first passage in~\eqref{pr5} ensues from~\eqref{hyp-k} and the
second one from the H\"older and Young inequalities, $\nu$ being a
positive constant to be chosen small enough. Furthermore,  also
taking into account~\eqref{form-b}, we have
\begin{equation}
\label{pr6}
\begin{aligned}
\int_0^t &  \tanh(r) b \left(\partial_{tt}^2 \ume \cdr, \partial_t
\ume \cdr \right)\dd r \\ &  \geq \frac{C_b}{2} \tanh(t) \|
\partial_t \ume (t) \|_{\bfw}^2 - \frac12 \int_0^t \tanh'(r)
b\left(\partial_{t} \ume \cdr, \partial_t \ume \cdr \right)\dd r \\
& \geq \frac{C_b}{2} \tanh(t) \| \partial_t \ume (t) \|_{\bfw}^2 -
\frac{K_b}{2}K_6^2\,,
\end{aligned}
\end{equation}
the last passage  due to~\eqref{aprio-long3-epsmu}
and~\eqref{form-b}. In the same way, we find
\begin{align}
&  \label{pr7} \int_0^t \tanh(r) \partial_{tt}^2\chime \cdr
\partial_{t}\chime \cdr \dd r \geq \frac{\tanh(t)}{2}
 \| \partial_t \chime (t) \|_{L^2 (\Gamma_c)}^2 - \frac{K_6^2}2\,,
 \\
 & \label{pr8}
\begin{aligned}
\int_0^t\int_{\Gamma_c} \tanh(r) \nabla \left(\partial_t \chime \cdr
\right) \nabla \left(\partial_t \chime \cdr \right)\dd r =\int_0^t
\tanh(r) \int_{\Gamma_c} |\partial_t (\nabla \chime\cdr )|^2 \dd
r\,,
\end{aligned}
\end{align}
and we remark that
\begin{equation}
\label{pr9}
\begin{aligned}
 & \int_0^t
\tanh(r) \langle \alpha'_\mu (\ume\cdr)
\partial_t \ume\cdr ,
\partial_t \ume\cdr\rangle \dd r \geq 0\,, \\ &  \int_0^t \tanh(r)
\int_{\Gamma_c} \beta_\mu'(\chime \cdr) |\partial_t \chime\cdr|^2
\dd r\geq 0\,.
\end{aligned}
\end{equation}
Finally, we perform the same computations as
in~\eqref{3.6.2.bis}--\eqref{3.6.3}, and observe that
\begin{align}
 &
\label{pr10}
\begin{aligned}
\int_{0}^t \int_{\Gamma_c} \tanh(r) \partial_t \left(\ppmu (\chime
\cdr) \ume \cdr \right) \partial_t \ume \cdr \dd r = I_{12}+
I_{13}\,,
\end{aligned}
\\
& \label{pr11}
 I_{12}= \int_{0}^t
\int_{\Gamma_c} \tanh(r)\ppmu (\chime \cdr) |
\partial_t \ume \cdr|^2  \dd r \geq 0\,,
\\
&  \label{pr12} \begin{aligned}
 I_{13} & = \int_{0}^t \tanh(r) \int_{\Gamma_c} \HHmu (\chime \cdr)
 \partial_t \chime \cdr \ume \cdr
\partial_t \ume \cdr \dd r \\ & \geq -  C \int_{0}^t \| \partial_t \chime
\|_{L^2 (\Gamma_c)} \| \ume  \|_{\bfw} \| \partial_t \ume \|_{\bfw}
\geq- C K_6^3\,,
\end{aligned}
\end{align}
the passages in~\eqref{pr12} due to the H\"older inequality, the
continuous embedding~\eqref{continuous-embedding}, and
estimates~\eqref{aprio-long3-epsmu}--\eqref{aprio-long4-epsmu}. In
the end, we point out that the analogues
of~\eqref{e:poincare}--\eqref{e:l4} and~\eqref{3.6.1}--\eqref{3.6.2}
go through. Then, collecting~\eqref{pr0}--\eqref{pr12} and putting
forth the same arguments as in the proof of
Proposition~\ref{prop:2.1}, we conclude estimate~\eqref{est:1-me}
for all $\delta>0$. We also find
\[
\| \partial_t \Thetame^{1/2}\|_{L^2 (\delta,+\infty;H)} + \|
\partial_t \Thetasme^{1/2}\|_{L^2 (\delta,+\infty;L^2 (\Gamma_c))} \leq
C_\delta\,,
\]
and, arguing in the same way as
in~\eqref{e:further-bound-1}--\eqref{e:holder}, we infer that for
all $\rho \in (0,2)$
\[
\| \partial_t \Thetame\|_{L^2 (\delta,+\infty;L^{12/7}(\Omega))} +
\|
\partial_t \Thetasme\|_{L^2 (\delta,+\infty;L^{2-\rho} (\Gamma_c))} \leq
C_{\delta,\rho}\,,
\]
whence the bound~\eqref{est:4-me} ensues.
 \fin
\subsection{Proof of Theorem~\ref{th:1.1}}
\label{ss:a.3}
\subsubsection{Passage to the limit in Problem $\pepsmu$ as $\eps \searrow 0$}
First, we introduce the boundary value problem obtained by taking
$\eps=0$ in Problem~$\pepsmu$, which we shall supplement with the
initial data~\eqref{cond-teta-zero}--\eqref{cond-chi-zero}.
\medskip\par\noindent
{\bf Problem $\pmu$}.
 Given a quadruple of initial data $(\teta_0,\teta_s^0,\uu_0,\chi_0)$
complying with~\eqref{cond-teta-zero}--\eqref{cond-chi-zero},
 find  functions  $( \teta, \teta_s,
\uu,\chi)$, $\uu$ with the regularity~\eqref{reguI}, $\chi$
with~\eqref{regchiI},  and, for all $T>0$,
\begin{align}
& \teta \in L^2 (0,T; V) \cap L^\infty (0,T;L^1(\Omega)), \ \
\calelleveps(\teta) \in H^1 (0,T;V')\,, \label{reg-teta-forte-bis}
\\
 & \teta_s \in   L^{2}(0,T;H^1 (\Gamma_c))  \cap L^{\infty}(0,T;L^1
(\Gamma_c)), \ \ \calelleveps(\teta_s) \in H^1 (0,T;H^1 (\Gamma_c)')
\,, \label{reg-teta-s-forte-bis}
\end{align}
 fulfilling the initial
conditions
\begin{equation}
\label{e:inimu}
\begin{aligned}
& \calelleveps (\teta(0)) = \calelleveps(\teta_0) \ \  \aein\,
\Omega\,, \qquad \calelleveps (\teta_s(0)) = \calelleveps(\teta_s^0)
\ \   \aein\, \Gamma_c\,,
\\
 & \uu(0)=\uu_0 \ \  \aein\,
\Omega\,, \qquad \chi(0)=\chi_0 \ \  \aein \,\Gamma_c\,,
\end{aligned}
\end{equation}
the equations
\begin{align}
& \label{teta-weak-app-bis}
\begin{aligned}
\langle
\partial_t\calelleveps(\teta) , v  \rangle -\int_{\Omega} \dive(\uu_t) \, v
+\int_{\Omega} \nabla \teta \, \nabla v  & + \int_{\Gamma_c} k(\chi)
(\teta-\teta_s) v\\ &  = \pairing{}{}{h}{v} \quad \forall\, v \in V
\ \hbox{ a.e. in }\, (0,+\infty)\,,
\end{aligned}
\\ & \label{teta-s-weak-app-bis}
\begin{aligned}
 \langle
\partial_t\calelleveps(\teta_s), v
\rangle  -\int_{\Gamma_c}
\partial_t \lambda(\chi) \, v  & +\int_{\Gamma_c} \nabla \teta_s  \, \nabla
v \\ &  = \int_{\Gamma_c} k(\chi) (\teta-\teta_s) v  \quad \forall\,
v \in \hunoc  \ \hbox{ a.e. in }\, (0,+\infty)\,,
\end{aligned}
\end{align}
and such that $(\uu,\chi)$ comply
with~\eqref{eqIappr}--\eqref{eqIIappr}.
\medskip
\\
Notice that no uniqueness result is available for (the Cauchy
problem for) Problem $\pmu$.

The following result is a straightforward consequence of
Proposition~\ref{prop:large-time}.
\begin{proposition}
\label{p:pass-1} Assume~\eqref{A5}--\eqref{hyp-k} and
 \eqref{e:consequences}--\eqref{cond-h-infty}.  Let $\mu>0$ be
 fixed. Then,
  there exists
 a (not relabeled) subsequence  of $\{
 (\tetame,\tetasme,\ume,\chime)\}_{\eps}$
 and  a quadruple $(\teta,\teta_s,
 \uu,\chi)$ such that  for all $T>0$ the following convergences
 hold as $\eps \searrow 0$
 \begin{subequations}
 \label{oh}
\begin{align}
\label{convtetae}
 &\tetame\weakto\teta\quad\text{in }L^2(0,T;V), \qquad
 \eps\mathcal{R}(\tetame)
 \to 0\quad\text{in }L^\infty(0,T;V'),
 \\
 &
\label{convtetase} \tetasme\weakto\teta_s\quad\text{in }L^2(0,T;H^1
(\Gamma_c)),  \qquad \eps\mathcal{R}_{\Gamma_c}(\tetasme)
 \to 0\quad\text{in }L^\infty(0,T;\hunoc'),
\\
& \label{conveweweak}
 \calelleveps(\tetame) \weakto \calelleveps(\teta)\quad \text{in }L^2(0,T;V)\,,
 \\
 & \label{convewestrong}
 \calelleveps(\tetame) \to \calelleveps(\teta)\ \  \text{in }\mathrm{C}^0([0,T];X) \  \ \text{for every Banach space
 $X$ with $V' \Subset X$}\,,
\\ &
\label{convewe} \eps \mathcal{R}(\tetame) +  \calelleveps(\tetame)
\weakto \calelleveps(\teta)  \quad \text{in }H^1(0,T;V')\,,
\\
& \label{convezeweak}
 \calelleveps(\tetasme) \weakto \calelleveps(\teta_s) \quad \text{in }L^2(0,T;\hunoc)\,,
 \\
 & \label{convezestrong}
 \calelleveps(\tetasme) \to \calelleveps(\teta_s) \ \  \text{in }\mathrm{C}^0([0,T];Y) \   \ \text{for every Banach space
 $Y$ with $\hunoc' \Subset Y$}\,,
\\ &
\label{conveze} \eps \mathcal{R}_{\Gamma_c}(\tetasme) +
\calelleveps(\tetasme) \weakto  \calelleveps(\teta_s)\quad \text{in
}H^1(0,T;\hunoc')\,,
\\
& \label{convchie}
\begin{aligned} &\chime\weaksto\chi\quad\text{in
}L^2(0,T;H^2(\Gamma_c)) \cap L^\infty(0,T;\hunoc)\cap
H^1(0,T;L^2(\Gamma_c)),
\\
&\chime\rightarrow\chi\quad\text{in } L^2(0,T;H^{2-\rho}(\Gamma_c))
\cap \mathrm{C}^0([0,T];H^{1-\rho}(\Gamma_c))\quad \text{for all
$\rho \in (0,1)$},
\end{aligned}
\\
& \label{convpichie}
\begin{aligned} &\ppmu(\chime)\weaksto\ppmu(\chi)\quad\text{in
} L^\infty(0,T;\hunoc)\cap H^1(0,T;L^2(\Gamma_c)),
\\
&\ppmu(\chime)\rightarrow\ppmu(\chi)\quad\text{in }
\mathrm{C}^0([0,T];H^{1-\rho}(\Gamma_c))\quad \text{for all $\rho
\in (0,1)$},
\end{aligned}
\\
 \label{convxie}
&\beta_\mu (\chime)\weakto\beta_\mu (\chi)\quad\text{in }L^2(0,T;L^2(\Gamma_c)),\\
 \label{convHe}
&\HHmu(\chime)\weaksto\HHmu(\chi)\quad\text{in
}L^\infty(0,T;L^p(\Gamma_c))  \quad \text{for all $ 1 \leq p < \infty$},\\
\label{convalphae}
&\alpha_\mu (\ume)\weakto \alpha_\mu (\uu)\quad\text{in }L^2(0,T;\bfw'),\\
\label{convue} &
\begin{aligned}
&
\ume\weakto\uu\quad\text{in }H^1(0,T;\bfw)\\
&\ume\to\uu\quad\text{in }\mathrm{C}^0([0,T];H^{1-\rho}(\Omega)^3) \quad
\text{for all $\rho \in (0,1)$\,,}
\end{aligned}
\end{align}\end{subequations}
(the symbol $\Subset$ in~\eqref{convewestrong}
and~\eqref{convezestrong} signifying compact inclusion),  and
$(\teta,\teta_s,\uu,\chi)$ is a solution to Problem~$\pmu$.

 Furthermore, $(\teta,\teta_s,\uu,\chi)$  has the additional
regularity~\eqref{e:further-reg1}--\eqref{e:further-reg3}
and~\eqref{e:further-reg4}
 on eve\-ry interval
$(\delta,T)$, for all $0<\delta<T$ and, up to the extraction of
another subsequence, for $\{
 (\tetame,\tetasme,\ume,\chime)\}_{\eps}$  the   enhanced
 convergences hold as $\eps \searrow 0$ for all $0<\delta <T$
 \begin{subequations}
 \label{ooh}
\begin{align}
\label{convtetame-en}
 &
\begin{aligned}
& \tetame \weaksto\,\teta\ \  \text{in }L^\infty(\delta,T;V)  \cap
H^1 (\delta,T; L^{12/7} (\Omega)), \\ &
 \tetame
\to \teta \ \ \text{in }\mathrm{C}^0 ([\delta,T];H^{1-\rho}(\Omega))
\ \text{for all }\rho \in (0,1),
\end{aligned}
 \\
 &
\label{convtetasme-en}
\begin{aligned}
 & \tetasme\weaksto\,\teta_s\quad\text{in
}L^\infty(\delta,T;H^1 (\Gamma_c)) \cap H^1 (\delta,T; L^{2-\rho}
(\Gamma_c))\
  \text{for all $\rho \in (0,2)$}\, \\ &
 \tetasme
\to \teta_s \ \ \text{in }\mathrm{C}^0
([\delta,T];H^{1-\rho}(\Gamma_c)) \ \text{for all }\rho \in (0,1),
\end{aligned}
\\
& \label{convume-en}
\begin{aligned}
 \ume\weaksto\uu\quad\text{in }W^{1,\infty}(\delta,T;\bfw)\,,
\end{aligned}
\\
& \label{convchime-en}
\begin{aligned} &\chime\weaksto\,\chi\quad\text{in
}L^\infty(\delta,T;H^2 (\Gamma_c)) \cap H^1(\delta,T;\hunoc)\cap
W^{1,\infty}(\delta,T;L^2(\Gamma_c)),
\\
&\chime\rightarrow\chi\quad\text{in
}\mathrm{C}^0([\delta,T];H^{2-\rho}(\Gamma_c))  \quad \text{for all
$\rho \in (0,2)$},
\end{aligned}
\\
 \label{convxien-en}
&\beta_\mu (\chime)\weaksto\,\beta_\mu (\chi)\quad\text{in
}L^\infty(\delta,T;L^2(\Gamma_c))\,.
\end{align}
\end{subequations}
\end{proposition}
\paragraph{Sketch of the proof.}
Since the proof of this result substantially goes along the same
lines as the proof of~\cite[Prop.~4.4]{bbr3}, we shall just outline
its main steps, referring to~\cite{bbr3} for all details.
\begin{itemize}
\item First of all, using estimates~\eqref{aprio-long-epsmu} and
standard compactness results (cf.~\cite{Simon87}), with a
diagonalization process we extract a subsequence of $\{
 (\tetame,\tetasme,\ume,\chime)\}_{\eps}$ for which
 convergences~\eqref{convtetae}--\eqref{convue} hold as $\eps \searrow
 0$. In order to show
that  $\teta \in L^\infty (0,T;L^1
 (\Omega))$,  as in the proof of~\cite[Thm.~1]{bbr3} we exploit a
 \emph{Lebesgue point argument}. Indeed, using the first
 of~\eqref{convtetae} and the convexity  of  $|\cdot|$,
  we
 find that for all $t_0 \in (0,T)$ and $r>0$
such that $(t_0-r,t_0+r) \subset (0,T)$
\begin{equation}
\label{e:lebesgue-point} \int_{t_0-r}^{t_0+r}\int_{\Omega} |\teta|
\leq \liminf_{\eps \searrow 0}
\left(\int_{t_0-r}^{t_0+r}\int_{\Omega} |\tetame| \right)\leq 2r
\sup_{\eps>0}\| \tetame \|_{L^\infty (0,T;L^1(\Omega))} \leq 2r
K_6\,,
\end{equation}
the latter inequality due to~\eqref{aprio-long2-epsmu}. We now
divide the above relation by $r$ and let $r \downarrow 0$. Using
that the Lebesgue point property  holds at almost every  $t_0 \in
(0,T)$, we
 obtain for all $T>0$ the estimate
\begin{equation}
\label{e:same-constant} \| \teta\|_{L^\infty (0,T;L^1(\Omega))} \leq
K_6\,.
\end{equation}
A completely analogous argument may be developed for proving that
$\teta_s \in L^\infty (0,T;L^1
 (\Gamma_c))$.
\item
 The identification of the limits of  the sequences
 $\{\beta_\mu(\chime)
 \}$ and $\{ \HHmu (\chime) \}$ follows from the strong-weak closedness
 of the graphs of the operators
 $\beta_\mu$, and $\HHmu$. As for the limits of  $\{ \alpha_\mu (\ume)\}$,
 $\{\calelleveps (\tetame)
 \}$ and
$\{\calelleveps (\tetasme)
 \}$,  thanks to~\cite[Lemma~1.3, p.~42]{barbu},
 it is  sufficient to prove that for all $t>0$
\[
\begin{aligned}
& \limsup_{\eps \searrow 0} \int_0^t \int_{\Omega} \calelleveps
(\tetame) \tetame \leq  \int_0^t \int_{\Omega} \calelleveps (\teta)
\teta\,,\\ &\limsup_{\eps \searrow 0} \int_0^t \int_{\Gamma_c}
\calelleveps (\tetasme) \tetasme \leq \int_0^t \int_{\Omega}
\calelleveps (\teta_s) \teta_s\,,
\\
&\limsup_{\eps \searrow 0} \int_0^t \int_{\Gamma_c} \alpha_\mu
(\ume) \cdot\ume \leq  \int_0^t \int_{\Gamma_c} \alpha_\mu (\uu)
\cdot\uu\,.
\end{aligned}
\]
Convergences~\eqref{convpichie} for $\{ \ppmu (\chime) \}$ ensue
from the corresponding~\eqref{convchie} for $\{ \chime \}$ and fron
the fact that $0 \leq \ppmu'(r) = \HH_\mu (r) \leq 1$ for all $r \in
\R$, cf.~\eqref{def-hmu}.
 Notice that the pair $(\uu,\chi)$
complies with the initial conditions~\eqref{iniu}--\eqref{inichi} in
view of convergences~\eqref{convue} and~\eqref{convchie}. In the
same way,
 combining
convergences~\eqref{tetazeroapp-4} and~\eqref{tetazeroapps-4} of the
sequences  $\{\calelleveps( \tetazeroemu)\}$, $\{
\calelleveps(\tetaessezeroemu)\}$, with
convergences~\eqref{convewestrong} and~\eqref{convezestrong},
respectively, we conclude the initial conditions~\eqref{e:inimu} for
$\calelleveps( \teta)$ and $\calelleveps( \teta_s)$. Finally, we
refer to the proof
 of~\cite[Prop.~4.4]{bbr3} for the argument showing that $(\teta,\teta_s,\uu,\chi)$ is a solution to
 Problem~$\pmu$.
\item
 The further
 regularity~\eqref{e:further-reg1}--\eqref{e:further-reg3} and~\eqref{e:further-reg4} for
 $(\teta,\teta_s,\uu,\chi)$
 and the enhanced convergences~\eqref{ooh} ensue from
 estimates~\eqref{est-delta-me},  see also the proof of
 Theorem~\ref{th:2.1}.
 \end{itemize}
\fin
\\
\noindent In view of the definition  of \emph{approximable solution}
to Problem~$\PP$ which we shall give in Section~\ref{ss:a.4},
 we select, among solutions of Problem~$\pmu$, only the
ones just constructed by passing to the limit in Problem~$\pepsmu$
as the viscosity parameter $\eps $ vanishes. Albeit improperly,
within the scope of this section we shall refer to them as
\emph{viscosity solutions}.
\begin{definition}[Viscosity solutions of Problem~$\pmu$]
\label{def:visc} Let $\mu>0$ be fixed. We say that a  quadruple
$(\teta, \teta_s, \uu,\chi)$ is a \emph{viscosity solution} of
Problem~$\pmu$ if
\begin{enumerate}
\item it is a solution to Problem~$\pmu$;
\item
 there
exist a sequence $\eps_k \searrow 0$ and a family of solutions
$\{(\tetamek,\tetasmek,\umek,\chimek) \}_k$ of Problem~$\pepskmu$
converging  as  $\eps_k \searrow 0 $ to $(\teta, \teta_s, \uu,\chi)$
in the sense specified by~\eqref{convtetae}--\eqref{convue}, on
every interval $(0,T)$.
\end{enumerate}
\end{definition}
\noindent In particular, it follows from the last part of
Proposition~\ref{p:pass-1} that every \emph{viscosity solution} to
Problem~$\pmu$ has the further
regularity~\eqref{e:further-reg1}--\eqref{e:further-reg3}
and~\eqref{e:further-reg4}.
\subsubsection{Passage to the limit as $\mu\searrow 0$}
\noindent
\begin{proposition}[Estimates on viscosity solutions of
Problem~$\pmu$] \label{prop:est-viscous}
Assume~\eqref{A5}--\eqref{hyp-k} and
 \eqref{e:consequences}--\eqref{cond-h-infty}.
Then,
\begin{enumerate}
\item
 estimates~\eqref{aprio-long-epsmu} on  the half-line $(0,+\infty)$, \eqref{est:1-me} on $(\delta,+\infty)$
  for
 all
$\delta>0$,  and~\eqref{est:4-me} on $(\delta,+\infty)$, for all
$\delta>0$ and $\rho \in (0,2)$,
 hold, with the same
constants $K_6$, $K_7 (\delta)$, and $K_8 (\delta,\rho)$,  for all
$\mu>0$ and for every viscosity  solution
$\{(\tetamu,\tetasmu,\uumu,\chimu) \}$ to Problem $\pmu$.
\item
Furthermore, for all $T>0$ there exists a constant $K_9 (T)$, only
depending on $T$, on  the quantity $M$~\eqref{def:M}, on the
functions $\lambda,$ $k$, and $\sigma$, as well as on $\|
\ln(\teta_0)\|_{H}$ and $\| \ln(\teta_s^0)\|_{L^2 (\Gamma_c)}$, such
that for all $\mu
>0$
\begin{equation}
\label{e:sorpresa} \| \calelleveps(\tetamu) \|_{L^\infty (0,T;H)} +
\| \calelleveps(\tetasmu) \|_{L^\infty (0,T;L^2(\Gamma_c))} \leq K_9
(T)\,.
\end{equation}
\end{enumerate}
\end{proposition}
\noindent As already mentioned in Remark~\ref{small-explanation},
the calculations leading to~\eqref{e:sorpresa} cannot be performed
on the solutions of Problem~$\pepsmu$, with $\eps>0$.
\paragraph{Sketch of the proof.} \, The first part of the statement
is an obvious consequence of the Definition~\eqref{def:visc} of
viscosity solution, of estimates~\eqref{aprio-long-epsmu}
and~\eqref{est-delta-me}, and of a trivial lower-semicontinuity
argument. The enhanced estimate~\eqref{e:sorpresa} for
$\calelleveps(\tetamu)$ and $\calelleveps(\tetasmu) $ is obtained by
testing~\eqref{teta-weak-app-bis} by $\calelleveps (\tetamu)$,
\eqref{teta-s-weak-app-bis} by $\calelleveps (\tetasmu)$, adding the
resulting relations and integrating on time. We refer
to~\cite[Sec.~4.2]{bbr3} for all details. Here, we just point out
that, being $\ln(\teta_0) \in H$,  thanks to~\eqref{e:inimu} there
holds
\[
\| \calelleveps (\tetamu (0))\|_{H} = \| \calelleveps
(\teta_0)\|_{H} \leq \| \ln(\teta_0)\|_H\,,
\]
and the analogous bound for $\| \calelleveps (\tetasmu (0))\|_{L^2
(\Gamma_c)} $ ensues from the fact that $\ln(\teta_s^0) \in L^2
(\Gamma_c)$. Such estimates are then used in the calculations
yielding~\eqref{e:sorpresa}.  \fin
\\
Thus, we have the following proposition, which concludes the proof
of Theorem~\ref{th:1.1}.
\begin{proposition}
\label{p:pass-2} Assume~\eqref{A5}--\eqref{hyp-k} and
 \eqref{e:consequences}--\eqref{cond-h-infty}.  Let
 $\{ (\tetamu,\tetasmu,\uumu,\chimu)\}_\mu$ be a sequence of \emph{viscosity solutions} to Problem~$\pmu$.
  Then,
  there exists
 a (not relabeled) subsequence  of $\{ (\tetamu,\tetasmu,\uumu,\chimu)\}_\mu$
 and  functions $(\teta,\teta_s,
 \uu,\chi,\eeta,\xi,\zeta)$ such that  for all $T>0$ the following convergences
 hold as $\mu \searrow 0$
 \begin{subequations}
 \label{ohhh}
\begin{align}
\label{convtetamu}
 &\tetamu\weakto\teta\quad\text{in }L^2(0,T;V),
 \\
 &
\label{convtetasmu} \tetasmu\weakto\teta_s\quad\text{in }L^2(0,T;H^1
(\Gamma_c)),
\\
& \label{conveweweamuk}
 \calelleveps(\tetamu) \weaksto \ln(\teta)\quad \text{in }L^\infty(0,T;H)\cap H^1(0,T;V') \,,
 \\
 & \label{convewemustrong}
 \calelleveps(\tetamu) \to \ln(\teta)\ \  \text{in }\mathrm{C}^0([0,T];X) \  \ \text{for every Banach space
 $X$ with $H \Subset X$}\,,
\\
& \label{convezemuweak}
 \calelleveps(\tetasmu) \weaksto \ln(\teta_s) \quad \text{in }L^\infty(0,T;L^2 (\Gamma_c))\cap H^1(0,T;\hunoc') \,,
 \\
 & \label{convezestrongmu}
 \calelleveps(\tetasmu) \to \ln(\teta_s) \ \ \text{in }\mathrm{C}^0([0,T];Y) \    \ \text{for every Banach space
 $Y$ with $L^2 (\Gamma_c) \Subset Y$}\,,
\\
& \label{convchiemu}
\begin{aligned} &\chimu\weaksto\chi\quad\text{in
}L^2(0,T;H^2(\Gamma_c))\cap L^\infty(0,T;\hunoc)\cap
H^1(0,T;L^2(\Gamma_c)),
\\
&\chimu\rightarrow\chi\quad\text{in }L^2(0,T;H^{2-\rho}(\Gamma_c))
\cap \mathrm{C}^0([0,T];H^{1-\rho}(\Gamma_c))   \quad \text{for all
$\rho \in (0,1)$},
\end{aligned}
\\
& \label{convpichiemu}
\begin{aligned} &\ppmu(\chimu)\weaksto\chi^+\quad\text{in
} L^\infty(0,T;\hunoc)\cap H^1(0,T;L^2(\Gamma_c)),
\\
&\ppmu(\chimu)\rightarrow\chi^+\quad\text{in }
\mathrm{C}^0([0,T];H^{1-\rho}(\Gamma_c))   \quad \text{for all $\rho
\in (0,1)$},
\end{aligned}
\\
 \label{convxiemu}
&\beta_\mu (\chimu)\weakto\xi \quad\text{in }L^2(0,T;L^2(\Gamma_c)),
 \qquad \xi \in \beta(\chi) \ \ \aein\, \Gamma_c \times (0,T),\\
 \label{convHemu}
&
\begin{aligned}
& \HHmu(\chimu)\weaksto\zeta\quad\text{in
}L^\infty(0,T;L^p(\Gamma_c)) \quad \text{for all $ 1 \leq p <
\infty$},
 \\ & \zeta \in \HH(\chi) \ \ \aein\, \Gamma_c \times (0,T),
\end{aligned}
 \\
\label{convalphaemu} &\alpha_\mu (\uumu)\weakto \eeta\quad\text{in
}L^2(0,T;\bfw'),
\qquad \eeta \in  \alpha(\uu) \ \aein\,  (0,T),\\
\label{convuemu} &
\begin{aligned}
&
\uumu\weakto\uu\quad\text{in }H^1(0,T;\bfw)\\
&\uumu\to\uu\quad\text{in }\mathrm{C}^0([0,T];H^{1-\rho}(\Omega)^3) \quad
\text{for all $\rho \in (0,1)$\,,}
\end{aligned}
\end{align}
\end{subequations}
and $(\teta,\teta_s,\uu,\chi,\eeta,\xi,\zeta)$ is a solution to
Problem~$\PP$.

 Furthermore, the functions $\teta,\,\teta_s,\,\uu, \,\chi$, and $\xi$  have the additional
regularity~\eqref{e:further-reg}
 on eve\-ry interval
$(\delta,T)$, for all $0<\delta<T$ and, up to the extraction of
another subsequence, for $\{
 (\tetamu,\tetasmu,\uumu,\chimu)\}_{\mu}$  the   enhanced
 convergences~\eqref{ooh} hold  as $\mu \searrow 0$ for all $0<\delta <T$.
\end{proposition}
\paragraph{Sketch of the proof.} The proof follows the very same lines of
the argument for Proposition~\ref{p:pass-1} (cf. also the proof
of~\cite[Thm.~1]{bbr3}). We just point out that
convergences~\eqref{convpichiemu} ensue from the pointwise
convergence $\ppmu(\chimu) \to \chi^+$ a.e. in $\Gamma_c \times
(0,T)$ (which can be deduced from~\eqref{convchiemu}, exploiting
formula~\eqref{formpmu}), and again from the bound~\eqref{def-hmu}
on $\ppmu'$. Furthermore,
\eqref{conveweweamuk}--\eqref{convezestrongmu} are a consequence of
the additional estimate~\eqref{e:sorpresa}. Finally, as in the proof
of Proposition~\ref{p:pass-1}, the strong
convergences~\eqref{convewemustrong} and~\eqref{convezestrongmu}, as
well as the Yosida convergences $\calelleveps(\teta_0) \weakto
\ln(\teta_0) $ in $H$ and $\calelleveps(\teta_s^0) \weakto
\ln(\teta_s^0) $ in $L^2 (\Gamma_c)$ as $\mu \searrow 0$,
 enable us to pass to the limit in initial
conditions~\eqref{e:inimu} and deduce that
\[
\ln(\teta(0))= \ln(\teta_0) \ \ \aein\ \Omega, \qquad
\ln(\teta_s(0))= \ln(\teta_s^0) \ \ \aein\ \Gamma_c\,,
\]
whence initial conditions~\eqref{iniw} and~\eqref{iniz}. \fin
\begin{remark}
\label{rem:precisazione} \upshape
 Notice that, in the statement of
Proposition~\eqref{p:pass-2},
assumptions~\eqref{cond-f-infty}--\eqref{cond-h-infty} on the data
$\mathbf{f}$, $g$, and $h$  are  stronger than the data requirements
in Theorem~\ref{th:1.1}, and there is the additional
condition~\eqref{e:consequences} on $\widehat{\beta}$. This is due
to the fact that, to avoid unnecessary repetitions,  we have chosen
to unify in the present section the proof of global existence with
the proof of the long-time estimates.

In fact,  a closer perusal of the proof of Theorem~\ref{th:1.1}  and
a comparison with the argument for~\cite[Thm.~1]{bbr3} show that the
sole~\eqref{A5}  on $\beta$ and~\eqref{hypo-h}--\eqref{hypo-g} on
the data $h$, $\mathbf{f}$, and $g$  are sufficient for   the global
existence of solutions to Problem~$\PP$ on the finite-time interval
$(0,T)$. Further, under  assumptions~\eqref{data-further-reg} (which
are the finite-time versions
of~\eqref{cond-f-infty}--\eqref{cond-h-infty}) one proves the
further regularity~\eqref{e:further-reg} on the interval $(0,T)$. In
particular, condition~\eqref{e:consequences} on $\widehat{\beta}$ is
not necessary for the aforementioned finite-time results.
\end{remark}

\subsection{Rigorous proof of  Proposition~\ref{prop:2.1}}
\label{ss:a.4} \noindent
 We are now in the position of specifying  approximable
solutions to Problem~$\PP$ as the ones arising as limits of
viscosity solutions to Problem~$\pmu$.
\begin{definition}[Approximable solutions of Problem~$\PP$]
\label{def:a2} We say that $(\teta, \teta_s,
\uu,\chi,\eeta,\xi,\zeta)$ is an \emph{approximable solution} of
Problem~$\PP$ if
\begin{enumerate}
\item  it is a
solution to Problem~$\PP$,
\item
 there
exist a sequence $\mu_k \searrow 0$ and  a family of viscosity
solutions $\{(\tetamuk,\tetasmuk,\uumuk,\chimuk) \}_k$ of
Problem~$\pmuk$ converging,  as $\mu_k \searrow 0 $,   to $(\teta,
\teta_s, \uu,\chi,\eeta,\xi,\zeta)$, in the sense specified
by~\eqref{convtetamu}--\eqref{convuemu}, on every interval $(0,T)$.
\end{enumerate}
\end{definition}
\begin{remark}
\label{rem:basta} \upshape Taking into account
Remark~\ref{rem:precisazione}, we might state Theorem~\ref{th:1.1}
in the following more precise way: under
assumptions~\eqref{A5}--\eqref{hypo-g}, on every interval $(0,T)$
Problem~$\PP$ admits at least an approximable solution $(\teta,
\teta_s, \uu,\chi,\eeta,\xi,\zeta)$. If, in addition,
\eqref{data-further-reg} holds, then every approximable solution
 has the further regularity~\eqref{e:further-reg}.
\end{remark}
\paragraph{Rigorous  proof of Proposition~\ref{prop:2.1}.}
It follows from the Definition~\eqref{def:a2} of approximable
solution, from the estimates on viscosity solutions to Problem
$\pmu$ specified in Proposition~\ref{prop:est-viscous}, and from
elementary  lower-semicontinuity arguments.
 \fin


                                %



\begin{thebibliography}{99}



\bibitem{barbu}
V.\ Barbu.
\newblock {\em Nonlinear Semigroups and Differential Equations in Banach
  Spaces}.
\newblock Noordhoff, Leyden, 1976.


\bibitem{bbr1}
E.~Bonetti, G.~Bonfanti, and R.~Rossi. \newblock Global existence
for a contact problem with adhesion.
\newblock {\em Math. Meth. Appl. Sci.}, 31, 1029--1064, 2008.



\bibitem{bbr2}
E. Bonetti, G. Bonfanti, and R. Rossi. \newblock Well-posedness and
long-time behaviour for a  model of contact with adhesion.
\newblock {\em
Indiana Univ. Math. J.},  56, 2787--2819, 2007.

\bibitem{bbr3} E. Bonetti, G. Bonfanti, and R. Rossi. \newblock
Thermal effects in adhesive contact: modelling and analysis.
\newblock Quaderno del Seminario Matematico di Brescia n. 07/2008, 2008.



\bibitem{bcfg1}
E. Bonetti, P. Colli, M. Fabrizio, and G. Gilardi. \newblock Global
solution to a singular integrodifferential system related to the
entropy balance. \newblock {\em Nonlinear Anal.},  66,  1949--1979,
2007.


\bibitem{bcfg2} E. Bonetti, P. Colli, M. Fabrizio, and G. Gilardi.
\newblock
Modelling and long-time behaviour for phase transitions with entropy
balance and thermal memory conductivity.
\newblock
{\em Discrete Contin. Dyn. Syst. Ser. B}, 6, 1001--1026, 2006.

\bibitem{bcf}
E. Bonetti, P. Colli, and M. Fr\'emond.
\newblock A phase field model with thermal memory governed by the entropy balance.
  \newblock {\em Math. Models Methods Appl. Sci.},   13, 1565--1588, 2003.


\bibitem{BFR}
E. Bonetti, M. Fr\'emond, and E.~Rocca.
\newblock A new dual approach for a class of phase transitions with
memory: existence and long-time behaviour of solutions.
\newblock {\em J. Math. Pures Appl.},
88,  455--481, 2007.


\bibitem{brezis73}
H.~Br\'ezis.
\newblock {\em Op\'erateurs Maximaux Monotones et Semi-groupes de Contractions
  dans les Espaces de Hilbert}.
\newblock Number~5 in North Holland Math. Studies. North-Holland, Amsterdam,
  1973.


\bibitem{colli92} P.~Colli.
\newblock On some doubly nonlinear evolution equations in Banach
spaces. \newblock \emph{Japan J. Indust. Appl. Math.},   9,
181--203, 1992.

\bibitem{colli-gilardi-laurencot-novickcohen}
P.~Colli, G.~Gilardi, P.~Lauren\c{c}ot, and A.~Novick-Cohen.
\newblock Existence and long-time behavior of the conserved
phase-field system with memory.
\newblock{\em Discrete Contin. Dynamic Systems}, 5, 375--390, 1995.



\bibitem{fs05}
E.~Feireisl and G.~Schimperna. \newblock Large time behaviour of
solutions to Penrose-Fife phase change models.
\newblock {\em Math. Methods Appl. Sci.}, 28, 2117--2132, 2005.



\bibitem{fre}
M.~Fr\'emond.
\newblock {\em Non-smooth Thermomechanics}.
\newblock Springer-Verlag, Berlin, 2002.


\bibitem{gps06}
M.~Grasselli, H.~Petzeltov\'{a}, and G.~Schimperna.
\newblock Long-time behavior of solutions to the Caginalp system
with singular potential. \newblock{\em Z. Anal. Anwend.}, 25,
51--72, 2006.


\bibitem{Haraux91} A.~Haraux. \newblock
Syst\`emes dynamiques dissipatifs et applications.
\newblock
RMA Res. Notes Appl. Math. 17. \newblock Masson, Paris, 1991.

\bibitem{Krejci-Zheng} P.~Krej{\v{c}}{\'{\i}} and  S.~Zheng.
\newblock Pointwise asymptotic convergence of solutions
 for a phase separation model.
\newblock
 {\em Discrete Contin. Dyn. Syst.},   16, 1--18,   2006.

 \bibitem{nuovaref}
P.~Rybka and K-H. Hoffmann,
\newblock Convergence of solutions to Cahn-Hilliard equation.
\newblock {\em Comm. Partial Differential equations}, 24, 1055--1077, 1999.


\bibitem{Simon87}
J.~Simon.
\newblock Compact sets in the space {$L^p(0,T; B)$}.
\newblock {\em Ann. Mat. Pura Appl. (4)}, 146, 65--96, 1987.




\end{thebibliography}
\end{document}